\documentclass[a4paper,12pt]{amsart}
\usepackage[utf8]{inputenc}
\usepackage[T1]{fontenc}
\usepackage[UKenglish]{babel}
\usepackage[margin=20mm]{geometry}
\usepackage{color}
\usepackage{bbm,bm}

\usepackage{subfigure}

\usepackage{graphicx}

\usepackage{amsmath,amssymb,amsfonts,amsthm}
\usepackage{mathrsfs,eucal,dsfont}
\usepackage{verbatim,enumitem}

\usepackage{hyperref,url}
\usepackage{amssymb, amsmath, amsthm, amscd}

\usepackage{eucal,mathrsfs,dsfont}
\usepackage{color}

\renewcommand{\leq}{\le}
\renewcommand{\geq}{\ge}

\newcommand{\e}{\text{e}}

\newcommand{\I}{\mathds 1}

\def\d{{\rm d}}
\def\<{\langle}
\def\>{\rangle}

\newtheorem{theorem}{Theorem}[section]
\newtheorem{lemma}[theorem]{Lemma}
\newtheorem{proposition}[theorem]{Proposition}

\numberwithin{equation}{section}

\theoremstyle{definition}
\newtheorem{definition}[theorem]{Definition}

\newtheorem{remark}[theorem]{Remark}
\newtheorem{assumption}{Assumption}

\begin{document}
\allowdisplaybreaks
\title[Heat kernel and exit time fluctuations on random graphs]
{\bfseries  Fluctuations of heat kernels and exit times 
on random graphs and application to 
the critical long-range percolation}

\author{Zherui Fan  \qquad  Lu-Jing Huang \qquad Takashi Kumagai}


\thanks{\emph{Z. Fan:}
School of Mathematical Sciences, Peking University, Beijing, China.
 \texttt{zheruifan@pku.edu.cn}}

\thanks{\emph{L.-J. Huang:}
School of Mathematics and Statistics \& Key Laboratory of Analytical Mathematics and Applications (Ministry of Education), Fujian Normal University, Fuzhou, China.
  \texttt{huanglj@fjnu.edu.cn}}

\thanks{\emph{T. Kumagai:}
Department of Mathematics, Waseda University, Shinjuku, Tokyo, 169-8555,
 Japan.\\
\texttt{t-kumagai@waseda.jp}}

\date{}
\maketitle

\begin{abstract}
Let $\mathcal{G}$ be a graph equipped with random conductance. 
We give sufficient conditions for the Markov chain on $\mathcal{G}$ to exhibit large scale fluctuations in its on-diagonal heat kernel and 
in exit times. The conditions also imply non-tightness of the height process w.r.t.\ the quenched measure, 
while it 
is tight w.r.t.\ the annealed measure.
As an application of the general theory, we prove large scale fluctuations for
the simple random walk on $1$-dimensional critical long-range percolation. 


\medskip


\end{abstract}
\allowdisplaybreaks


\section{Introduction}
Over the past several decades, random geometry has made tremendous progress.
Detailed geometric properties of random media, especially at \lq criticality\rq,\, 
have been investigated. 
As research in random geometry has advanced, questions about the dynamics on these media have also attracted attention. In particular, the asymptotic behavior of random walks on (discrete) random media has been studied. It turns out that these walks often behave anomalously, with fluctuations arising in many forms.

Let us begin with a fundamental example: the incipient infinite cluster (IIC) for critical percolation on Galton-Watson trees. 
For simplicity, we consider the case of a binomial offspring distribution, which can be viewed as 
bond percolation on a binary tree.
Let $\mathbb B$ be the binary tree rooted at $0$. 
Let $\{\eta_e: e \mbox{ is a bond on } \mathbb B\}$ be 
{\it i.i.d.} random variables with $P[\eta_e=1]=\frac 12$, $P[\eta_e=0]=\frac 12$ for each $e$; this defines critical bond percolation on $\mathbb B$. Set
$$ {\mathcal C}(0) = \{x \in \mathbb B: \hbox{ there is an $\eta$-open path between $0$ and $x$}\}. $$
Let $d(\cdot,\cdot)$ denote graph distance on $\mathbb B$, and set 
\[\mathbb B_n=\{x\in \mathbb B: d(0,x)=n\},\qquad 
\mathbb B_{\le n}=\{x\in \mathbb B: d(0,x)\le n\}.\] 
It is easy to see that 
$W_n := | {\mathcal C}(0) \cap \mathbb B_n|$ is a Galton-Watson branching process
with offspring distribution ${\rm Bin}(2,\frac 12)$. 
Since $E[W_1]=1$, 
the process ${W_n}$ becomes extinct almost surely, and hence $|\mathcal C(0)|<\infty$ almost surely. 
However, in any ball of radius $n$, there are open clusters with diameter of order $n$ with high probability. To study mesoscopic properties of these large, finite clusters, we regard them as subsets of an infinite cluster defined as follows. First, we construct the limiting probability measure of the process conditioned on non-extinction. Namely, for $A\subset \mathbb B_{\le k}$,
$$\lim_{n\to\infty} P\left[{\mathcal C}(0)\cap \mathbb B_{\le k} =A\mid W_n \ne 0\right]
 = |A \cap \mathbb B_k|\, P\left[{\mathcal C}(0)\cap \mathbb B_{\le k} =A\right]
 =:\mathds{P}_{\rm IIC}(A).$$
By a standard extension theorem, $\mathds P_{\rm IIC}$ extends uniquely to a probability measure $\mathds P$ on the set of infinite connected subsets of $\mathbb B$ containing $0$. This measure is called the incipient infinite cluster (IIC).

Let $\mathcal G=\mathcal G(\omega)$ be the rooted, labeled tree with distribution $\mathds P$. In 1986, Kesten \cite{Kes86a} studied simple random walk
$(\{X_t\}_{t\ge 0}, \mathbf{P}_x=\mathbf{P}_x^\omega, x\in {\mathcal G})$ 
on $\mathcal G$ and proved that the walk behaves anomalously. (More precisely, \cite{Kes86a} considered discrete-time random walk.) In particular, it was shown that, as $n\to\infty$,
$\{Z_t^{(n)}:=t^{-1/3}d(0,X_{nt})\}_{t\ge 0}$ converges weakly under $\mathds P\times \mathbf{P}_0$, 
where $d(\cdot,\cdot)$ denotes graph distance on $\mathcal G$. Barlow and Kumagai \cite{BK06} established more detailed results on the behavior of $\{X_t\}_{t\ge 0}$. They obtained detailed estimates for the random walk’s heat kernel $p_t^\omega(\cdot,\cdot)$; in particular, they proved the following fluctuation of the heat kernel:
\begin{equation}\label{eq: logloginf}
 \liminf_{t \to \infty} 
   (\log\log t)^{1/6} t^{2/3} p^\omega_{2t}(0,0)  \le 2, 
\quad \mathds P\text{-a.s.}\end{equation}
They also obtained annealed heat-kernel estimates. Specifically, for $t\ge c(d(x,y)\wedge 1)$, the following holds \cite[Theorem 1.5]{BK06}:
$$
c_1t^{-2/3}
\exp\left\{-c_2\left(d(x,y)^3/t\right)^{1/2}\right\}
\le
\mathds E_{x,y}\big[p_t^\omega(x,y)\big]
\le
c_3t^{-2/3}
\exp\left\{-c_4\left(d(x,y)^3/t\right)^{1/2}\right\},
$$
where ${\mathds E}_{x,y}[\cdot]:={\mathds P}[\cdot\mid {x,y}\in \mathcal G]$.
Thus, \eqref{eq: logloginf} indeed reflects fluctuations of the heat kernel. It was also shown that the processes $\{Z_t^{(n)}\}_{t\ge 0}$ defined above are not tight under the quenched measure $\mathbf{P}_0=\mathbf{P}_0^\omega$.
In fact, \eqref{eq: logloginf} follows readily from the following volume fluctuation result \cite[Proposition 2.8(a)]{BK06}: for every $\varepsilon>0$,
$$  \limsup_{n \to \infty} \frac{ V(0,n)}{ n^2(\log\log n)^{1-\varepsilon}}=\infty, 
 \quad \mathds P\text{-a.s.},$$
where $\bm{V}(0,n):=\sharp \{z\in  {\mathcal G}: d(0,z)\le n\}$. 
In other words, information about the random geometry yields information about the dynamics on the medium.
Similar heat kernel fluctuations have been obtained for the IIC of critical infinite-variance branching processes \cite{CroKum}, and for $2$-dimensional (resp. $3$-dimensional) uniform spanning trees in \cite{BCK21} (resp. \cite{DS-SW}). Intuitively, volume fluctuations arise because the scaling limit of the random medium (or random graph) $\mathcal G$ is itself random. Indeed, in the examples above, the scaling limits are random media. One important implication of such fluctuations is that the family of rescaled random walks is not tight (see Proposition \ref{thm:tightnot}). However, establishing these fluctuations in any specific example is generally highly non-trivial.

We note that heat kernel fluctuations of this kind do not occur in non-critical cases. As a typical example, consider supercritical bond percolation on $\mathds Z^d$, where $d\ge 2$, and let $\mathcal C_\infty(\omega)$ be the unique infinite open cluster. Barlow \cite{Bar04} proved the following quenched heat kernel estimates: there exist constants $c_1,\ldots,c_4>0$ and a finite random variable $S_x(\omega)$ for each $x\in\mathds Z^d$ such that, almost surely,
\begin{equation}\label{eq:HKGauss}
c_1t^{-d/2}\exp\Big\{-c_2|x-y|^2/t\Big\}\le p_t^\omega(x,y) \le 
c_3t^{-d/2}\exp\Big\{-c_4|x-y|^2/t\Big\}, 
\end{equation}
for $x,y\in {\mathcal C}_\infty (\omega)$ and $t\ge S_x(\omega)\vee |x-y|$. 
Similar heat kernel estimates have been established for more general random conductance models on 
$\mathds Z^d$ (see \cite{ABDH, ADS} etc.). 
\eqref{eq:HKGauss} implies that $c_1\le t^{d/2}p_t^\omega(x,y)\le c_3$.
Moreover, a quenched local central limit theorem has been proved \cite{BH09}, so fluctuations such as those in \eqref{eq: logloginf} do not occur. Note that, in this example, the scaling limit of $\mathcal C_\infty(\omega)$ is $\mathds R^d$, which is deterministic.

Although heat kernel and exit time fluctuations are important properties linking random geometry to the analysis of random walks on random graphs, no general theory has yet been presented in a clear and unified form. 
(We note that for processes on continuum random media, there is a nice work by Croydon \cite{Croydon07}.) 
The aim of this paper is to identify geometric sufficient conditions for these fluctuations and apply the resulting theory to critical long-range percolation (LRP).
In the next section, we first develop a general theory of heat kernel and exit time fluctuations, and then apply it to $1$-dimensional critical LRP. As we discuss in Section \ref{LRP-subsec}, critical LRP is a very challenging model, and its geometric and analytic properties have only recently begun to be investigated in detail; see \cite{Baumler23a, Baum23b, DFH23,DFH25+,FH26}. The main tool for analyzing the medium is the so-called multiscale analysis. In this paper, we substantially extend this method and obtain sharp estimates for volume and resistance. 
In particular, the main technical novelties are as follows.
\begin{itemize}
\item[(1)] We sharpen resistance tail estimates by controlling the dependence of moment bounds on the moment order and iteratively improving exponential moment estimates. This yields exponential upper-tail bounds for resistance diameters with any exponent below $1/(1-\delta)$, where $\delta$ is the resistance scaling exponent; see Lemma~\ref{improve-eta} and Proposition~\ref{uppertail-up}. Conversely, by suppressing long edges and forcing escaping flows through successive blocks, we obtain a probability lower bound for large escape resistance with exponent $1/(1-\delta)$; see Proposition~\ref{uppertail-0}. These estimates give nearly matching tail exponents despite the presence of long-range shortcuts.

\item[(2)] We construct intervals consisting of $2K$ blocks of length $m$, in which the resistance from the origin to every vertex is of order at most $m^\delta$. The key is a hierarchical arrangement of long edges with nested attachment intervals; see Lemma~\ref{lem-sepcont}. By adding buffer intervals and restricting certain long edges, we also ensure a large resistance from the central interval to the exterior. These properties hold simultaneously with probability at least $\exp\{-CK\}$; see Lemma~\ref{lem-singlesc}. These estimates form the basis of the multiscale construction of resistance balls with large volume and large escape resistance.

\item[(3)] We handle the dependence caused by long edges through multiscale arguments. For large-volume events, we obtain independent trials by contracting previously explored regions and control the errors from contraction and long-range connections to transfer the estimates to the original infinite network; see Lemma~\ref{lem-probAk}. For small-volume events, the main difficulty is to combine resistance conditions with opposite monotonicities in the edge configuration. We overcome this using conditional probability estimates, the FKG inequality and the improved resistance diameter tail bounds; see Subsection~\ref{subsect-proofThm}  for more details. These arguments establish almost-sure volume fluctuations and quantitative control of the scales of large-volume events; see Propositions~\ref{Vn-lu} and~\ref{prop-assumpC2}.
\end{itemize} 
These estimates should be useful in future studies of the model and are also of independent interest.

There are many examples of random walks in random media at criticality. To understand their long-time behavior, it would be useful to establish fluctuations in their on-diagonal heat kernels and exit times. In particular, we leave the following as an open problem concerning random walks on IICs; see \cite{Kum14} for definitions and known results on the IICs and random walks:

\medskip
\begin{center}
{\bf Open problem}: Prove Assumption \ref{assumpB} and \ref{assumpC} in Section \ref{general-subsec} 
for simple random walk on the IIC of bond percolation on ${\mathds Z^d}$ for $d\ge 11$ and for $d=2$. 
\end{center}
\medskip
The remainder of the paper is organized as follows. In Section 2.1, we develop the general theory of heat kernel and exit time fluctuations. In Section 2.2, we review the history of research on LRP and state our main theorem on fluctuations for $1$-dimensional critical LRP. Section 3 is devoted to the proof of Theorem \ref{osc-LRP}, the main result on heat kernel and exit time fluctuations for random walks on $1$-dimensional critical LRP. Some of the detailed estimates in this section may also be useful in future research on the model. Finally, Section 4 proves the results on the general theory of fluctuations stated in Section 2.1, using methods and techniques from the theory of resistance forms.

\section{Framework and main results}\label{main-results}
\subsection{General theory}\label{general-subsec}
Consider a connected undirected graph $\mathcal{G}=(V,E)$ with vertex set $V$ and edge set $E$. 
We assume that $\mathcal{G}$ is locally finite, meaning that $\text{deg}(x):=\#\{y\in V:\ \langle y,x\rangle\in E \}<\infty$ for all $x\in V$. 
We also assume that there is a distinguished point $0\in V$. 
For convenience, we will also use $\sim$ to denote edges, specifically, $x\sim y$ implies $\langle x,y\rangle\in E$. 
The edges are endowed with positive conductances $\mu_e=\mu_{xy}$ for $e=\langle x,y\rangle \in E$. 
Denote $\mu_x=\sum_{y\sim x}\mu_{xy}$ for each $x\in V$ and set 
\begin{equation}\label{def-meas}
\mu(A)=\sum_{x\in A}\mu_x\quad  \text{for each }A\subset V.
\end{equation}

We define a quadratic form $(\mathcal{D},\mathds{H}^2)$ as follows:
\begin{equation*}
  \begin{aligned}
  \mathcal{D}(f,g)=\frac{1}{2}\sum_{x\sim y}\mu_{xy}(f(x)-f(y))(g(x)-g(y)), ~~~
  \mathds{H}^2=\left\{f\in \mathds{R}^{V}:\ \mathcal{D}(f,f)<\infty\right\}. 
  \end{aligned}
\end{equation*} 
Let $X=(X_t)_{t\geq 0}$ be the continuous-time random walk on $\mathcal{G}$ associated with $\mathcal{D}$ and the measure $\mu$. This means that $X$ is the Markov process with generator 
\begin{equation}\label{def-generator}
  \mathcal{L}f(x):=\frac{1}{\mu_x}\sum_{y\sim x}\mu_{xy}(f(y)-f(x)). 
\end{equation}
Denote $\mathbf{P}_x$ and $\mathbf{E}_x$ as the probability and expectation of $X$ starting from $x$. We also define the transition density (heat kernel density) of $X$ with respect to $\mu$ as 
\begin{equation}\label{def-hk}
   p_t(x,y)=\mathbf{P}_x[X_t=y]/\mu_y. 
\end{equation}
It is easy to check that $p_t(\cdot,\cdot)$ is symmetric, that is, $p_t(x,y)=p_t(y,x)$ for all $x,y\in V$ and all $t>0$.

For any two disjoint subsets $A,B\subset V$, we define the effective resistance between $A$ and $B$ as 
\begin{equation}\label{def-eff}
    1/R(A,B)=\inf\left\{\mathcal{D}(f,f):\ f\in \mathds{H}^2,\ f|_A=1\text{ and }f|_B=0\right\}. 
\end{equation}
This formula is known as the Dirichlet principle for effective resistance. 
Additionally, there exists a dual form of the Dirichlet principle referred to as Thompson's principle of effective resistance; see \cite{DS84}. 
To introduce this, we begin with some notations from electrical networks. 
Let $f$ be a function defined on $V\times V$. We call $f$ a flow on the graph $\mathcal{G}$ if it satisfies 
\begin{equation*}
  f_{xy}=
  \begin{cases}
    -f_{yx}, & \mbox{if } \langle x,y\rangle \in E,\\
    0, & \mbox{otherwise}
  \end{cases}
\end{equation*}
for all $x,y\in V$. Thus, the net flow out of a vertex $x$ is given by $f_{(x)}=\sum_{y\neq x}f_{xy}$. 
Moreover, given two disjoint subsets $A,B\subset V$, we say $f$ is a unit flow from $A$ and $B$ if it is a flow satisfying 
\begin{equation*}
  \sum_{x\in A}f_{(x)}=1\quad \text{and}\quad f_{(y)}=0\quad \text{for all }y\notin A\cup B. 
\end{equation*}
Then the effective resistance $R(A,B)$ can also be expressed as 
\begin{equation}\label{def-eff-flow}
   R(A,B)=
  \inf\left\{\frac{1}{2}\sum_{x\sim y} f_{xy}^2/\mu_{xy}:\ f \mbox{ is a unit flow from $A$ to }B\right\}.
\end{equation}
In particular, when $A=\{x\}$ is a singleton, we write $R(x,B)$ for $R(\{x\},B)$. 
It is important to note that the effective resistance $R(\cdot,\cdot)$ defines a metric on the graph $\mathcal{G}$; see \cite[Section 2.3]{Ki01} for more details. 
For any $x\in V$ and $r>0$, let $B(x,r)$ denote the ball centered at $x$ with radius $r$ under the resistance distance, that is,
\begin{equation}\label{def-ball}
  B(x,r)=\left\{y\in V:\ R(x,y)\leq r\right\}.
\end{equation}
We also denote $\bm{V}(x,r)=\mu(B(x,r))$ as the volume of the ball $B(x,r)$ with respect to $\mu$.

In the following, consider a probability space $(\Omega, \mathcal{F},\mathds{P})$ carrying a family of random graphs $\mathcal{G}(\omega)=(V(\omega), E(\omega)),\ \omega\in \Omega$. 
We assume that for each $\omega\in \Omega$, the graph $\mathcal{G}(\omega)$ is locally finite, connected, and includes a marked vertex $0\in V(\omega)$. 
Let $R(\cdot,\cdot)=R(\cdot,\cdot)(\omega)$ denote the effective resistance on $\mathcal{G}(\omega)$. We denote the resistance balls in $\mathcal{G}(\omega)$ by $B(x,r)=B(x,r)(\omega)$ and their volume by $\bm V(x,r)=\bm V(x,r)(\omega)$. 
We then assign positive conductances $\mu_{e}=\mu_{e}(\omega)=\mu_{xy}(\omega)$ for each edge $e\in E(\omega)$, thereby forming a random electric network denoted as $(\mathcal{G},\{\mu_{e}\}_{e\in E})$. 

For each $\omega\in \Omega$, we let $X=((X_t)_{t\geq 0}, (\mathbf{P}_x^\omega)_{x\in V(\omega)})$ represent the continuous-time random walk on $\mathcal{G}(\omega)$ with generator given by  \eqref{def-generator}, and let $p_t^\omega(\cdot,\cdot)$ denote the transition density of $X$ with respect to $\mu$, as defined in \eqref{def-hk}.

We now introduce the following assumptions regarding the vertices in the graph $\mathcal{G}$, which will help us in establishing the tightness of the on-diagonal heat kernels.

\begin{assumption}\label{assumpA}
We say a vertex $x\in V$ satisfies Assumption A if there exists a constant $d_f>0$ and a decreasing function $h:\mathds{R}_+\rightarrow \mathds{R}_+$ with $\lim_{\lambda\rightarrow \infty}h(\lambda)=0$ such that the following holds. 
  \begin{itemize}
     \item[(A1)] There exist constants $q_1,q_2\geq 1$  such that  for all $r>0$ and $\lambda\geq 1$,  
  \begin{equation*}
    \mathds{P}\left[\lambda^{-q_1}r^{d_f}\leq \bm V(x,r)\leq \lambda^{q_2}r^{d_f}\right]\geq 1- h(\lambda). 
  \end{equation*}
  
  \item[(A2)] For all $r>0$ and $\lambda\geq 1$, 
  \begin{equation*}
 \mathds{P}\left[R(x,B(x,r)^c)\geq \frac{r}{2\lambda}\right]
\geq 1-h(\lambda). 
  \end{equation*}
  \end{itemize}
\end{assumption}

\begin{theorem}\label{thm-tight}
  Suppose that $x\in V$ satisfies Assumption \ref{assumpA}.  Then there exist constants $C_1<\infty$ and $\eta>0$ {\rm(}depending on $q_1$ and $q_2${\rm)} such that for all $\varepsilon \in (0,1)$ and all $r, t\geq 1$,
\begin{equation}\label{Exit-tight}
  \mathds{P}\left[\varepsilon \leq
  r^{-1-d_f}\mathbf{E}_x\left[\tau_{B(x,r)}\right]\leq
  \varepsilon^{-1}\right]\geq  1-C_1h(\varepsilon^{-\eta}),
\end{equation}
\begin{equation*}\label{hk-tight}
  \mathds{P}\left[\varepsilon \leq  t^{\frac{d_f}{1+d_f}}p_{2t}^\omega(x,x)\leq \varepsilon^{-1}\right]\geq  1-C_1h(\varepsilon^{-\eta}),
 \end{equation*}
where $\tau_A$ denotes the first exit time of the random walk $X$ from $A$.  
\end{theorem}

Next, we introduce conditions that will imply fluctuation in the heat kernel. To this end,  
let ${\bf SV}:=\{\phi:\ \mathds{N}\rightarrow [0,\infty) \,|\, \phi$ is strictly increasing 
with $\phi(1)=1$, satisfying \eqref{cond-phi}$\}$, where 
\begin{equation}\label{cond-phi}
\lim_{n\rightarrow \infty}\phi(n)=\infty\quad \text{and}\quad\lim_{n\rightarrow \infty}\frac{\phi(n)}{n^{\varepsilon}}=0
\end{equation}
for all $\varepsilon\in(0,1)$. 
For convenience, for each $\phi \in {\bf SV}$, we set $\phi(0)=0$, $\phi(\infty)=\infty$, extending it to $\phi:\ [0,\infty]\rightarrow [0,\infty]$ such that $\phi$ is continuous and strictly increasing. 

\begin{assumption}\label{assumpB} (1) We say the random graph $\mathcal{G}$ satisfies Assumption B(1) if 
there exist $c_1,d_f>0$, an increasing 
sequence $\{r_n\}_{n\geq 1}\subset {\mathds N}$ and a strictly increasing function $h_1\in  {\bf SV}$ such that 
  \begin{equation*}
 \mathds{P}[\{ E_n \text{ i.o.}\}]=1, ~~~\text{ where } E_n:= \left\{{\bm V(0,r_n)}\ge c_1r_n^{d_f} h_1(r_n)\right\}.
  \end{equation*} 
(2) We say the random graph satisfies Assumption B(2) if 
there exist $c_2,d_f>0$, an increasing 
sequence $\{r_n\}_{n\geq 1}\subset {\mathds N}$ and a strictly increasing function $h_2\in  {\bf SV}$ such that 
  \begin{equation*}
 \mathds{P}[\{ F_n \text{ i.o.}\}]=1, ~~~\text{ where } F_n:= \left\{{\bm V(0,r_n)}\le c_2r_n^{d_f} /h_2(r_n)\right\}.
  \end{equation*} 
\end{assumption}

Before presenting the main results under Assumption \ref{assumpB}, let us first provide a remark regarding this assumption.

 \begin{remark} Clearly, Assumption \ref{assumpB}(1) 
  implies 
   \begin{equation*}
    \limsup_{n\rightarrow\infty}\frac{\bm V(x,n)}{n^{d_f}h_1(n)}\ge c_1\quad \mathds{P}\text{-a.s.},\qquad \text{hence}\quad \limsup_{n\rightarrow\infty}\frac{\bm V(x,n)}{n^{d_f}h_1(n)^{1-\varepsilon}}=\infty\quad \mathds{P}\text{-a.s.}
  \end{equation*} 
  for any $\varepsilon\in (0,1)$. Similarly, Assumption \ref{assumpB}(2) 
  implies 
  \begin{equation}\label{eq:liminfvol}
    \liminf_{n\rightarrow\infty}\frac{\bm V(x,n)h_2(n)}{n^{d_f}}\le c_2\quad \mathds{P}\text{-a.s.},\qquad \text{hence}\quad      \liminf_{n\rightarrow\infty}\frac{\bm V(x,n)h_2(n)^{1-\varepsilon}}{n^{d_f}}
=0\quad \mathds{P}\text{-a.s.}
  \end{equation} 
  for any $\varepsilon\in (0,1)$.
\end{remark}

Under Assumption \ref{assumpB}, we have the following. 

\begin{theorem}\label{hk-limif}
{\rm(1)} Suppose that the random graph $\mathcal{G}$ satisfies Assumption \ref{assumpB}{\rm(1)}.  Then we have 
  \begin{equation*}
    \liminf_{t\rightarrow \infty}t^{\frac {d_f}{1+d_f}}\Phi_1(t)p_{2t}^\omega(0,0)\leq 2,\quad \mathds{P}\text{-a.s.}, 
  \end{equation*}
  where $\Phi_1\in {\bf SV}$ is determined by
  $h_1(x^{1/(1+d_f)}/\Phi_1(x))=\Phi_1(x)^{1+d_f}$.\\
 {\rm(2)}  Suppose that the random graph $\mathcal{G}$ satisfies Assumption \ref{assumpB}{\rm(2)}.  Then we have 
  \begin{equation*}
    \liminf_{r\rightarrow \infty}r^{-(1+d_f)}h_2(r)\mathbf{E}_0\left[\tau_{B(0,r)}\right]
    \leq 2c_2,\quad \mathds{P}\text{-a.s.} 
  \end{equation*}
   \end{theorem}

The corresponding $\limsup$ 
statements of the heat kernel and exit times require 
much stronger assumption. 
For this purpose, we first define \lq good sets\rq \  as follows.

\begin{definition}\label{def-I-1}
For $r,d_f>0$ and $\lambda, u_1,u_2\geq 1$, define the events
\begin{eqnarray*}
H_1(r,\lambda)&=&\left\{{R(0,B(0,r)^c)\geq r/{2\lambda}}\right\},\\
 H_2(r,u_1)&=&\left\{\bm V(0,r)\leq c_3r^{d_f}/u_1\right\},\\
 H_3(r,\lambda, u_2)&=&\left\{\bm V(0,\frac{r}{4\lambda})\geq c_4\Big(\frac{r}{{ 4}\lambda}\Big)^{d_f}/u_2\right\},
 \end{eqnarray*}
for some $c_3,c_4>0$.
\end{definition}

\begin{assumption}\label{assumpC}
(1) We say the random graph $\mathcal{G}$ satisfies Assumption C(1) if 
there exist an increasing sequence $\{r_n\}_{n\geq 1}\subset {\mathds N}$, 
a constant $d_f>0$ and strictly increasing functions $\lambda, h_2,h_3\in  {\bf SV}$ that satisfies
\begin{equation}\label{condC}
G(n):=\frac{h_2(n)^3}{\lambda(n)^{2+3d_f}h_3(n)^{2+\frac{d_f}{1+d_f}}}\in {\bf SV}
\end{equation}
such that 
  \begin{equation}\label{eq:iocond}
 \mathds{P}\left[\{ H_1(r_n,\lambda (n))\cap H_2(r_n,h_2(n))\cap H_3(r_n,\lambda (n),h_3(n))\text{ i.o.}\}\right]=1.
  \end{equation} 
  (2) We say the random graph $\mathcal{G}$ satisfies Assumption C(2)  if 
there exist an increasing sequence $\{r_n\}_{n\geq 1}\subset {\mathds N}$, 
a constant $d_f>0$ and strictly increasing functions $\lambda, h_3\in  {\bf SV}$ that satisfies
\begin{equation}\label{condD2}
\widehat G(n):=\frac{h_3(n)}{\lambda(n)^{1+d_f}}\in {\bf SV}
\end{equation}
such that 
  \begin{equation*}
 \mathds{P}\left[\{ H_1(r_n,\lambda (n))\cap   H_3(r_n,\lambda (n),1/h_3(n))\text{ i.o.}\}\right]=1.
  \end{equation*} 
\end{assumption}
Note that $E_n=H_3(4\lambda_nr_n,1/h_1(r_n))$ and $F_n=H_2(r_n,h_2(r_n))$, which are the events in Assumption \ref{assumpB}. 
In the following, we often write $\lambda_n:=\lambda (n)$. Let
\begin{equation}\label{eq:T_ndesu}
T(n)=T_n:=\Big(\frac{r_n}{\lambda_n}\Big)^{1+d_f}h_3(n)^{-1}.   
\end{equation}
Since $r_n\ge n$ and $\lambda ,h_3\in {\bf SV}$, $\lim_{n\to\infty}T(n)=\infty$. So we may 
assume that $\{T_n\}$ is an increasing sequence (otherwise, we take an increasing subsequence). 

\medskip 
Under this assumption, we have the following. 

\begin{theorem}\label{hk-limsup2}
{\rm(1)} Suppose that the random graph $\mathcal{G}$ satisfies Assumption \ref{assumpC}{\rm(1)}.  Then 
there exists a non-random constant $C_2>0$ such that  
  \begin{equation*}
    \limsup_{t\rightarrow \infty}\frac{t^{\frac {d_f}{1+d_f}}p_{2t}^\omega(0,0)}{\Phi_2(t)}\ge C_2,\qquad \mathds{P}\text{-a.s.}, 
  \end{equation*}
  where $\Phi_2\in {\bf SV}$ is determined by
  $\Phi_2(x)=G(T^{-1}(x))$  
  and $T^{-1}$ is the inverse function of $T$.  \\
  {\rm(2)} Suppose that the random graph $\mathcal{G}$ satisfies Assumption \ref{assumpC}{\rm(2)}.  Then 
there exists a non-random constant $C_3>0$ such that  
  \begin{equation*}
  \limsup_{r\rightarrow \infty}
  \frac{\mathbf{E}_0\left[\tau_{B(0,r)}\right]}{r^{1+d_f}\Phi_3(r)}
    \geq C_3,\qquad \mathds{P}\text{-a.s.} 
  \end{equation*}
  where $\Phi_3\in {\bf SV}$ is determined by
  $\Phi_3(x)=\widehat G(r^{-1}(x))$ 
  and $r^{-1}$ is the inverse function of $r$.  
  \end{theorem}
  
\medskip
Let
\[\widehat H_1(r,\lambda):=\left\{R(z,B(0,r)^c)\geq r/\lambda\quad \text{for all }z\in B(0,r/(2\lambda))\right.\}.
\]
Then $ H_1(r,\lambda)\subset \widehat H_1(r,\lambda)$.  Indeed, 
let $\bar R(\cdot,\cdot)$ be the resistance metric of the network that shorten
  $B(0,r)^c$. Noting that $\bar R(z,B(0,r)^c)=R(z,B(0,r)^c)$, we have 
  \begin{equation}\label{eq:resttri}
  R(z,B(0,r)^c)=\bar R(z,B(0,r)^c)\ge \bar R(0,B(0,r)^c)-\bar R(z,0)\ge R(0,B(0,r)^c)- R(z,0),
  \end{equation}
  where we use the shorting law of the resistance in the last inequality. In order to prove Assumption \ref{assumpC}, it is clearly enough to prove \eqref{eq:iocond} for $\widehat H_1(r,\lambda)$ instead of 
$H_1(r,\lambda)$, that is we will do in the rest of the paper. So, with some abuse of notation, we 
write $ H_1(r,\lambda)$ for $\widehat H_1(r,\lambda)$ below.

  \begin{remark} (1) Theorems \ref{hk-limif} and \ref{hk-limsup2} imply that  
  the order of fluctuations of the long-time behavior of $t^{\frac {d_f}{1+d_f}}p_{2t}^\omega(0,0)$ is 
  at least $\Phi_1\vee \Phi_2$. 
  Indeed, suppose $t^{\frac {d_f}{1+d_f}}\mathds{E}[p_{2t}^\omega(0,0)]\ge C$(resp.\ 
  $\le C$), then Theorem \ref{hk-limif} 
  (resp.\ \ref{hk-limsup2}) implies that the order of fluctuations is at least $\Phi_1$ (resp.\ $\Phi_2$).\\
  (2) The choices of $r_n,\lambda_n$ and $h_1(n), h_2(n), h_3(n)$ depend on models and
  they are crucial. For the random walk on long-range percolation with $d=1$ and $s=2$, 
  as we see in the following sections, 
  in order to check Assumption \ref{assumpB}(1) for instance, we choose
 $r_n=(\log n)^{-\alpha}M^{\delta n^2}$, $h_1(n)= (\log\log n)^{\alpha'}$ and to check Assumption \ref{assumpC}(1), we choose $r_n=M^{n^2}, \lambda_n=(\log n)^\iota$ and $h_2(n)= (\log n)^{\frac{(1-\varepsilon)(1-\delta)}{\delta}}$, 
  $h_3(n)=(\log n)^{1/(\delta\eta)}$, where $\delta (=1/d_f) \in (0,1)$ is the exponent of the growth of the effective resistance (see Theorem \ref{eff-tightness} below) in this model.  
  In the previous works, 
  Theorems \ref{hk-limif}(1) and \ref{hk-limsup2}(1) (for some case only Theorem \ref{hk-limif}(1)) are  obtained for random walks on the incipient infinite cluster (IIC) 
  of critical $n_0$-ary branching processes in \cite{BK06},
  on the IIC of critical infinite variance branching processes in \cite{CroKum}, and on  
  $2$-dimensional (resp.\ $3$-dimensional) uniform spanning trees in \cite{BCK21}
  (resp.\ \cite{DS-SW}). In these works, the graph distance is used, but since the models
  are homogeneous trees, the graph distance coincides with the resistance metric. \\
 (3) For deterministic and continuum media, results corresponding to 
Theorem \ref{hk-limif}(1) and Theorem \ref{hk-limsup2}(1) are established in  
\cite[Theorem 2]{Croydon07} (see also \cite[Theorem 3]{Croydon07}). \\
(4) In this paper, we use the resistance metric for the estimates. 
 It is possible to generalize the results that hold for some other metric (such as a graph distance) 
 under some additional conditions that relate the resistance metric and the metric in question; see for instance, \cite{Kum14,KM08}.\\   
(5) It is possible to state a version of Theorem \ref{hk-limsup2}(1) under the following assumption instead
of Assumption \ref{assumpC}(1);  
\begin{equation*}\label{condD1}
G(n):=\frac{h_3(n)^{2+\frac{d_f}{1+d_f}}}{\lambda(n)^{2+3d_f}h_2(n)^3}\in {\bf SV}\end{equation*}
such that 
  \begin{equation*}
 \mathds{P}\left[\{ H_1(r_n,\lambda (n))\cap H_2(r_n,1/h_2(n))\cap  H_3(r_n,\lambda (n),1/h_3(n))\text{ i.o.}\}\right]=1.
  \end{equation*} 
Since the proof is similar and we do not use it, we omit details. 
 \end{remark}
  
  \medskip
  
  Finally, we give a criterion for non-tightness of the process in the quenched sense. To this end, 
  define the rescaled  height process as 
  \[
  Z_t^{(n)}=n^{-\frac 1{d_f+1}}R(0, X_{nt}),\quad n\in \mathds{N},~t\ge 0. 
  \]
 By \eqref{Exit-tight}, $\{ Z_t^{(n)}\}_{n\in \mathbb N}$ 
are tight w.r.t.\ the annealed law 
given by  $\mathds{P}\times {\mathbf P}$. However, due to the large scale fluctuations, we do not have quenched tightness as in the following proposition. 
   \begin{proposition}\label{thm:tightnot}
  Under Assumption \ref{assumpB}{\rm(2)}, $\mathds{P}$-a.s., $\{Z^{(n)}: n\ge 1\}$ are not tight 
  w.r.t.\ $ \mathbf{P}_0$. 
  \end{proposition}

The proofs of  Theorems \ref{thm-tight}, \ref{hk-limif}, \ref{hk-limsup2}
and Proposition \ref{thm:tightnot} will be presented in Section \ref{sect-proofmr}. 
The primary technique for this proof involves establishing relationships between the volume of resistances and the heat kernel. We note that this method has been applied to the random walk on the incipient infinite cluster of trees and oriented percolation in \cite{BK06,BJKS08}. In this context, we use this approach to establish fluctuations for random walks on general random graphs. 

\subsection{Long-range percolation}\label{LRP-subsec}
The goal of this subsection is to apply the main results outlined in Section \ref{main-results} to the long-range percolation (LRP) model and to establish oscillatory behavior for the random walk within this model. 
Before that, we will provide a survey of the relevant research related to this model.  

Long-range percolation (LRP), initially introduced in \cite{S83,ZPL83}, is a percolation model on $\mathds{Z}^d$ in which each pair of vertices can be connected with a bond. This distinguishes LRP from classical Bernoulli percolation, as it allows for the existence of ``long-range jumps'', which leads to complex phase transition phenomena.
To be more precise, consider the sequence $\{p_x\}_{x\in \mathds{Z}^d}$, where $p_x\in [0,1]$ and $p_x=p_{x'}$ for all $x,x'\in \mathds{Z}^d$ such that $|x|=|x'|$. We assume that 
\begin{equation*}
  0<\beta:=\lim_{|x|\rightarrow \infty}p_x|x|^s<\infty 
\end{equation*}
for some $s>0$. 
The LRP model on $\mathds{Z}^d$ is defined by edges $\langle x,y\rangle$ occurring independently with probability $p_{x-y}$. 
In particular, $p_1$ represents the probability that two nearest neighbors are connected. 

The study of the LRP model has sparked great interests. 
One of the key topics is the case where $p_1<1$, which addresses the important question of the existence and uniqueness of the infinite component.
In one dimension, the first relevant work was conducted by \cite{S83}, which established that percolation does not occur if $s>2$.
\cite{NS86} demonstrated the existence of oriented percolation when $s<2$. 
Moreover, it was discovered that non-oriented percolation occurs in the critical case $s=2$ if $\beta$ is sufficiently large and $p_1$ is close to 1. This result for $s=2$ was further improved by \cite{AN86}, which established the criticality of $\beta$. Specifically, it showed that percolation does not occur if $\beta\leq 1$, whereas non-oriented percolation does occur if $\beta>1$ and $p_1$ is close to 1. 
Afterwards, \cite{MSV10} also proved that $\beta_c=1$ for the oriented case with $s=2$. 
In the high-dimensional case, the introduction of long-range connections is not crucial for the existence of percolation. Instead, the main focus shifts to quantifying the impact of these connections on the critical behavior; see e.g. \cite{MS96,SSV99}. 
Subsequently, people have also shown interest in the limiting properties of the graph distance metric when an infinite cluster is present. 
It is known that the renormalization structure leads to five distinct regimes, depending on whether $s$ is smaller, greater or equal to $d$ or $2d$. The graph distance exhibits different behaviors in each of these regimes; for a comprehensive overview of this topic, see e.g.\ \cite{DFH23}. 
Notably, in the critical case where $s=2d$, \cite{DFH23} proved that the subsequential limit of the graph distance is uniquely characterized by a natural list of axioms, based on the work in \cite{Baumler23a} regarding the tightness of the graph distance.

Another intriguing direction in the study of the LRP model is the behavior of the random walk on its infinite cluster.
Significant progress has been made in understanding the properties of random walks in the non-critical cases (i.e.\ $(d,s)\neq (d,2d)$) of the model. 
For instance, \cite{M08} provided some up-to-constant estimates for the box-to-box resistance when $d=1,\ s\in (1,2)\cup (2,\infty)$ as well as for $d\geq 2$ and $s\in (d,d+2)\cup (d+2,\infty)$. 
Additionally, \cite{KM08} employed methods to estimate volumes and effective resistances, ultimately deriving the corresponding heat kernel estimates for the case where $d=1,\ s>2$ and $p_1=1$. 
Furthermore, \cite{CS12} provided some bounds of the heat kernel and the spectral gap of the random walk for the case $d\geq 1,\ s\in (d,(d+2)\wedge 2d)$. 
Recently, \cite{CCK22} refined heat kernel estimates with bounds that allow for log-corrections across all non-critical scenarios. 
Moreover, the scaling limits of random walks on the non-critical LRP model were studied by \cite{CS13} and \cite{BT24}. 
Specifically, \cite{CS13} showed that for $s\in (d,d+1)$, the random walk on the infinite cluster of the LRP model converges to an  $(\alpha=s-d)$-stable process. 
Notably, in the case where $d=1$ and $s>2$, \cite{CS13} showed that the random walk converges to Brownian motion. These results were later generalized by \cite{BT24} to the case of $d\geq2$ and $s\in [d+1,d+2)$, showing that the associated random walk also converges to an $(\alpha=s-d)$-stable process.  
Let us emphasize that for these non-critical cases, there are no fluctuations of the heat kernel as in 
\eqref{eq: logloginf}, and the scaling limits exist.

Now, let us turn our attention to random walks on the critical (i.e.\ $(d,s)=(d,2d)$) LRP model. 
The first important result in this area was presented by \cite{B02}, 
which established recurrence for the random walk when $s=2d$ for $d=1,2$ by showing that the effective resistance from the origin to $[-N,N]^c$ diverges as $N\rightarrow \infty$. 
This result was subsequently extended in 
\cite{Baum23b} 
to more general LRP models with weight distributions satisfying certain moment conditions. 
In addition, \cite{M08} provided up-to-constant estimates for the box-to-box resistance in the critical LRP model when $d\geq 3$. Most recently, \cite{DFH24+,DFH25+} establish polynomial growth of the effective resistance for the critical LRP model in the case of $d=1$. 
Furthermore, \cite{CCK22} also provided precise heat kernel estimates that allow for log-corrections, as well as spectral dimensions for the random walks on the critical LRP models when $d\geq 2$.  
Very recently, \cite{FH26} extended this result to one-dimensional critical LRP model. 
Specifically, the authors presented non-Gaussian heat kernel estimates and spectral dimensions for the random walk on the one-dimensional critical LRP model, using tail estimates of the effective resistance from \cite{DFH25+}.

Despite these progresses, it is still a challenging problem to study the scaling limit of the random walks in the critical LRP models. 
As an application of our general results, we aim to explore this topic further. 
In particular, we intend to show that the heat kernels of the random walk on the one-dimensional critical LRP model exhibit oscillatory behavior, which implies that the random walk in this model does not have a scaling limit in the quenched sense.

To achieve this, we will focus on the one-dimensional critical LRP model in the following. Specifically, let edges $\langle i,j\rangle$ with $|i-j|=1$ (i.e. $i$ and $j$ are nearest neighbors) occur independently with probability $p_1=1$, and let edges $\langle i,j\rangle$ with $|i-j|>1$ occur independently with probability 
\begin{equation}\label{prob-pij}
  p_{i,j}=1-\exp\left\{-\beta\int_i^{i+1}\int_j^{j+1}\frac{1}{|u-v|^2}\d u\d v\right\}. 
\end{equation}
Here, $\beta>0$ is a parameter of this LRP model. We also call this model as a $\beta$-LRP model, and denote $\mathcal{E}$ for the edge set. 
Note that we chose the expression \eqref{prob-pij} for the connecting probability due to its strong scaling invariance property, making it a convenient choice. It is important to note that our proofs extend to cases where the connection probability is within a multiplicative constant of $\beta |i-j|^{-2}$.  

By placing a unit conductance on each edge $e\in \mathcal{E}$, i.e. $\mu_e\equiv 1$, we may then view this LRP model as a random electric network on $\mathds{Z}$. 
In this context, we have $\mu_x=\sum_{y\sim x}\mu_{xy}=\text{deg}(x)$ for $x\in \mathds{Z}$. 
Recall that $R(\cdot,\cdot)$ denotes the effective resistance as defined in \eqref{def-eff} and \eqref{def-eff-flow}. 
For any two disjoint subsets $I_1,I_2\subset \mathds{Z}$, we define $I_1\sim I_2$ as the event that $I_1$ and $I_2$ are directly connected by some long edges, and let $I_1\nsim I_2$ be the complement of this event, meaning that $I_1$ and $I_2$ are not directly connected by any long edge.

The following result is from \cite{DFH25+}, which establishes the polynomial growth of the effective resistance in one-dimensional critical LRP model. This result will 
play an important role in proving the oscillatory behavior of the heat kernel of the associated random walk in this model.

\begin{theorem}[{\cite[Theorem 1.1]{DFH25+}}]\label{eff-tightness}
  For all $\beta>0$, there exists an exponent $\delta=\delta(\beta)\in (0,1)$ {\rm(}depending only on $\beta${\rm)} such that for all $n\in \mathds{N}$, 
  \begin{equation*}
    \mathds{E}\left[R(0,[-n,n]^c)\right]\asymp_P R(0,[-n,n]^c)\asymp_P  n^\delta, 
  \end{equation*}
  and conditioned on $[-n,n]\nsim [-2n,2n]^c$, 
    \begin{equation*}
    \mathds{E}\left[R([-n,n],[-2n,2n]^c)\mid [-n,n]\nsim [-2n,2n]^c\right]\asymp_P R([-n,n],[-2n,2n]^c)\asymp_P  n^\delta, 
  \end{equation*}
  where $a(n)\asymp_P b(n)$ means that for all $\varepsilon>0$, there exist constants $0<c<C<\infty$ such that $\mathds{P}[cb(n)\leq a(n)\leq Cb(n)]\geq 1-\varepsilon$ for all $n\in \mathds{N}$. 
\end{theorem}

As an application of Theorems \ref{thm-tight}, \ref{hk-limif}, \ref{hk-limsup2} and \ref{eff-tightness}, we will show that the heat kernel of the random walk on the $\beta$-LRP model exhibits 
oscillatory behavior. 
The next theorem is our main theorem for the random walk on the $\beta$-LRP model. 
\begin{theorem}\label{osc-LRP}
For all $\beta>0$, the following holds. 
\begin{itemize}
  \item[\rm(1)] There exists a constant $q>0$ {\rm(}{depending only on $\beta$\rm)} such that for all $x\in \mathds{Z}$, 
   $r,t\ge 1$  
  and sufficiently small $\varepsilon \in(0,1)$, 
      \begin{eqnarray*}
 && \mathds{P}\left[\varepsilon \leq
  r^{-\frac{1+\delta}{\delta}}\mathbf{E}_x\left[\tau_{B(x,r)}\right]\leq
  \varepsilon^{-1}\right]\geq  1-\varepsilon^q,\\
  &&  \mathds{P}\left[\varepsilon \leq  t^{\frac{1}{1+\delta}}p_{2t}(x,x)\leq \varepsilon^{-1}\right]\geq  1-\varepsilon^q.  
      \end{eqnarray*}
  \item[\rm(2)]   There exist constants $\alpha_1,\alpha_2,C_4>0$ {\rm(}{depending only on $\beta$\rm)} such that for all $x\in \mathds{Z}$, 
    \begin{eqnarray*}
    \liminf_{t\rightarrow \infty}t^{\frac{1}{1+\delta}}(\log \log t)^{\alpha_1}p_{2t}^\omega(x,x)\leq 2,&& \mathds{P}\text{-a.s.},\\
 \liminf_{r\rightarrow \infty}r^{-\frac{1+\delta}\delta}(\log \log r)^{\alpha_2}\mathbf{E}_x\left[\tau_{B(x,r)}\right]\leq C_4, && \mathds{P}\text{-a.s.} 
  \end{eqnarray*}
  \item[\rm(3)] There exist constants $\alpha_3,\alpha_4>0$ {\rm(}{depending only on $\beta$\rm)} such that for all $x\in \mathds{Z}$,  
    \begin{eqnarray}
    \limsup_{t\rightarrow \infty} t^{\frac{1}{1+\delta}}(\log \log t)^{-\alpha_3}p_{2t}(x,x)=\infty, &&\mathds{P}\text{-a.s.}, \label{limsuppt}\\
 \limsup_{r\rightarrow \infty}r^{-\frac{1+\delta}\delta}(\log \log r)^{\alpha_4}\mathbf{E}_x\left[\tau_{B(x,r)}\right]=\infty, && \mathds{P}\text{-a.s.} \label{limsupext}
 \end{eqnarray}
  \item[\rm(4)] Let $X=(X_t)_{t\geq 0}$ be the continuous-time random walk on the $\beta$-LRP model and let $Z^{(n)}_t=n^{-\frac{\delta}{1+\delta}}R(0,X_{nt})$ for $n\geq 1$ and $t\geq 0$. 
 Then $\{Z^{(n)}:n\geq 1\}$ are tight w.r.t. the annealed measure 
 $\mathds{P}\times \mathbf{P}_0$, but $\mathds{P}$-a.s.\  they 
 are not tight w.r.t.\ the quenched measure $\mathbf{P}_0$. 
\end{itemize}
\end{theorem}

\section{Proof of Theorem \ref{osc-LRP}}

The objective of this section is to verify that the $\beta$-LRP model satisfies Assumptions \ref{assumpA}, \ref{assumpB} and \ref{assumpC}. Consequently, by the results in Section \ref{main-results}, we will show that the random walk on this model exhibits oscillatory behavior as presented in Theorem \ref{osc-LRP}.

Recall that $\mu_e\equiv 1$ for all $e\in \mathcal{E}$ in the $\beta$-LRP model.  
For $x\in \mathds{Z}$ and $r>0$, $B(x,r)$ denotes the ball centered at $x$ with radius under the resistance distance, and $\bm{V}(x,r)$ denotes the volume of the ball $B(x,r)$ with respect to $\mu$.  That is, 
\begin{equation*}
  \bm{V}(x,r)=\sum_{y\in B(x,r)} \mu_y=\sum_{y\in B(x,r)} \text{deg}(y). 
\end{equation*}
We also recall that $\delta=\delta(\beta)>0$ is the exponent defined in Theorem \ref{eff-tightness}. 
A main input of checking Assumptions \ref{assumpA}, \ref{assumpB} and \ref{assumpC}(1) is to show that fluctuations in $n^{-1/\delta}\bm{V}(0,n)$ as follows. 

\begin{proposition}\label{Vn-lu}
For all $\beta>0$, we have the following two statements.
\begin{itemize}

\item[\rm (1)] There exists a constant $\alpha_1>0$ {\rm(}depending only on $\beta${\rm)} such that 
  \begin{equation*}\label{volume-large}
    \limsup_{n\to\infty} \frac{\bm{V}(0,n)}{n^{1/\delta}(\log\log n)^{\alpha_1}}=\infty\quad \mathds{P}\text{-a.s.}
  \end{equation*}

\item[\rm (2)] There exist constants $\alpha_2, c_1>0$ {\rm(}both depending only on $\beta${\rm)} such that 
  \begin{equation*}\label{volume-small}
    \liminf_{n\to\infty}\frac{\bm{V}(0,n)(\log\log n)^{\alpha_2}}{n^{1/\delta}}\leq c_1 \quad \mathds{P}\text{-a.s.}
  \end{equation*}
  \end{itemize}
\end{proposition}

In addition, the following proposition shows that the $\beta$-LRP model satisfies Assumption 
\ref{assumpC}(2). 
To do this,  
fix $\eta>0$, $0<\theta<1$, and \begin{equation}\label{def-alpha}
\alpha>1/(1-\theta) 
\end{equation}
in the following. 
Let $g:\mathds N\to[1,\infty)$ be deterministic and nondecreasing, with
\begin{equation}\label{def-functg}
g(n)\rightarrow\infty
\quad\text{and}\quad 
g(n)^{1+d_{f}}=o(\log n), 
\end{equation}
where $d_f:=1/\delta$ and $\delta\in (0,1)$ is the exponent in Theorem \ref{eff-tightness}.

For $r\geq 0$ and $\lambda,h\geq 1$, let $\mathcal H(r,\lambda,h)$ denote the event that $H_1(r,\lambda)\cap H_3(r,\lambda,1/h)$, i.e., 
\begin{equation}\label{def-mathcalH}
\mathcal H(r,\lambda,h)
=\left\{R(0,B(0,r)^c)\geq \frac{r}{2\lambda},\  \bm V(0,\frac{r}{4\lambda})\geq h\left(\frac{r}{4\lambda}\right)^{d_f}\right\}
\end{equation}

\begin{proposition}\label{prop-assumpC2}
For all $\beta>0$, there exist constants $\lambda_0\geq2$ and $C_0<\infty$, depending only
on $\beta$, such that, with
\[
h(n):=
1+\frac{\theta}{C_0}\log\left(\lfloor n^\alpha\rfloor+{\rm e}\right)
\quad\text{and}\quad 
\lambda(n):=\lambda_0g(n),
\]
almost surely there exists a strictly increasing sequence $(r(n))_{n\geq1}$
such that
$
\mathcal H(r(n),\lambda(n),h(n))
$
holds for every $n\geq1$ and
\begin{equation}\label{lrp:target-speed}
\frac{\log r(n)}{n^\alpha\log n}
\rightarrow
\alpha(1+\eta)\delta\quad \text{as }n\rightarrow\infty.
\end{equation}
Moreover,
\begin{equation}\label{lrp:target-ratio}
\frac{h(n)}{\lambda(n)^{1+d_{f}}}
=
\frac{\alpha\theta}{C_0\lambda_0^{1+d_{f}}}
\frac{\log n}{g(n)^{1+d_{f}}}
\,(1+o(1))
\rightarrow\infty\quad \text{as }n\rightarrow\infty.
\end{equation}
\end{proposition}

The remainder of this section is organized as follows.
Subsections \ref{subsect-proof312} and \ref{better-up-tail} provide the resistance and degree estimates used below. 
In Subsection \ref{subsect-largeball}, we give the proof of Proposition \ref{prop-assumpC2}, from which Proposition \ref{Vn-lu}(1) follows.
Finally, in Subsection \ref{subsect-proofThm}, we verify Assumption \ref{assumpC}(1), deduce Proposition \ref{Vn-lu}(2), and complete the proof of Theorem \ref{osc-LRP}.

\subsection{Estimates for the resistance and degrees}\label{subsect-proof312}

In this subsection, 
we need some notation for effective resistance restricted to an interval.
Specifically,  
recall \eqref{def-eff-flow} that 
for any interval $I\subseteq \mathds{Z}$ and any two disjoint subsets $I_1,I_2\subseteq I$, 
$R_I(I_1,I_2)$ is 
the effective resistance between $I_1$ and $I_2$ restricted to the interval $I$, that is, 
\begin{equation*}\label{int-eff}
\begin{aligned}
   R_I(I_1,I_2)&:= 
  \inf\bigg\{\frac{1}{2}\sum_{x\sim y}f_{xy}^2:\ f\text{ is a unit flow between $I_1$ and $I_2$} \\
  &\qquad \qquad \qquad \qquad \qquad\text{and $f_{xy}=0$ for all }\langle x,y\rangle\in \mathcal{E}\setminus\mathcal{E}_{I\times I}\bigg\},
  \end{aligned}
\end{equation*}
where $\mathcal{E}_{I\times J}:=\mathcal{E}\cap (I\times J)$ for $I,J\subseteq \mathds{Z}$. In particular, when $I_1=\{i\}$ is a singleton, we also write $R_I(i,I_2)$ for $R_I(I_1,I_2)$.
Define the diameter of $I$ with respect to $R(\cdot,\cdot)$ by  
\begin{equation*}\label{def-diam}
\operatorname{diam}(I;R)=\max_{x,y\in I}R_I(x,y). 
\end{equation*}

We begin with the following upper tail of the resistance, which is a main input of verifying Assumption \ref{assumpC}(1).

\begin{proposition}\label{uppertail-0}
  For all $\beta>0$, there exists a constant $c_2>0$ {\rm(}depending only on $\beta${\rm)} such that for all $n\in \mathds{N}$ and $\varepsilon\in (0,1)$ with $\varepsilon^{1/(1-\delta)}n\geq 1$,
  \begin{equation*}
    \mathds{P}\left[R(0,[-n,n]^c)\geq \varepsilon^{-1}n^\delta\right]\geq \exp\left\{-c_2\varepsilon^{-\frac{1}{1-\delta}}\right\}. 
  \end{equation*}
\end{proposition}

To prove Proposition \ref{uppertail-0}, we need some decompositions of the resistance. To achieve this, fix $n\in \mathds{N}$ and $\varepsilon\in (0,1)$. Let $B_n$ be the event that $0\nsim [-n,n]^c$, and let $f$ be a unit flow from 0 to $[-n,n]^c$. Then, we can see that the flow $f$ emanates from 0 and ultimately flows into $[-n,n]^c$.
As a result, on the event $B_n$, the flow $f$ must pass through intervals $[-n,0)$ or $(0,n]$.
So let us define $\theta_1$ as the portion of flow $f$ that passes through the interval $(0,n]$ and then flows into $[-n,n]^c$. Similarly, we define $\theta_2$ by replacing $(0,n]$ with $[-n,0)$ in the definition of $\theta_1$.
Note that the flow may go back and forth between the intervals $[-n,0)$ and $(0,n]$. However, the definitions of $\theta_1$ and $\theta_2$ here only focus on the interval from which the flow exits $[-n,n]$. We also let
\begin{equation*}
  |\theta_1|=\sum_{i\in (0,n]}\sum_{j\in [-n,n]^c}f_{ij} \quad \text{and}\quad
  |\theta_2|=\sum_{i\in [-n,0)}\sum_{j\in [-n,n]^c}f_{ij}
  \end{equation*}
represent amounts of flows exiting from intervals $(0,n]$ and $[-n,0)$, respectively. It is clear that $\max\{|\theta_1|,|\theta_2|\}\geq 1/2$. In addition, define
\begin{equation}\label{def-Rhat}
  \begin{aligned}
  \widehat{R}([-n,0],[-n,n]^c)
   &:=
  \inf\bigg\{\frac{1}{2}\sum_{u\sim v}f_{uv}^2:\ f \text{ is a unit flow from $[-n,0]$ to $[-n,n]^c$ }, \\
  &\quad \quad \quad \text{and $f_{uv}=0$ for all }\langle u,v\rangle\in \mathcal{E}_{[-n,0]\times [-n,n]^c} \bigg\}
  \end{aligned}
\end{equation}
as the effective resistance generated by unit flows from $[-n,0]$ to $[-n,n]^c$, passing through the interval $(0,n]$.
Here, for any subsets $I_1,I_2\subset \mathds{Z}$, $\mathcal{E}_{I_1\times I_2}:=\mathcal{E}\cap(I_1\times I_2)$. 
 Similarly, denote $\widehat{R}([0,n],[-n,n]^c)$ as the effective resistance generated by 
unit flows from $[0,n]$ to $[-n,n]^c$, passing through the interval $[-n,0)$. Then we can see that
\begin{equation}\label{R-Rhat}
  R(0,[-n,n]^c)\geq \max\left\{|\theta_1|^2\widehat{R}([-n,0],\,[-n,n]^c),\,  |\theta_2|^2\widehat{R}([0,n],\,[-n,n]^c)\right\}.
\end{equation}

Let $B_{n,\varepsilon}$ be the event that there is no edge longer than $\varepsilon^{\frac{1}{1-\delta}}n$ (the red dashed curves in Figure \ref{event-B}) in the interval $[-n,n]$. 
Note that $B_{n,\varepsilon}\subset B_n$.
By the FKG inequality, we have that
\begin{equation}\label{prob-B}
\begin{aligned}
  \mathds{P}\left[B_{n,\varepsilon}\right]\geq \mathds{P}\left[0\nsim \left[-\varepsilon^{\frac{1}{1-\delta}}n,\varepsilon^{\frac{1}{1-\delta}}n\right]^c \right]^{2n+1}
  &=\exp\left\{-2(2n+1)\beta\int_0^1\int_{\varepsilon^{\frac{1}{1-\delta}}n}^{+\infty}\frac{1}{|x-y|^2}\d x\d y\right\}\\
  &\geq \exp\left\{-6\beta \varepsilon^{-\frac{1}{1-\delta}}\right\}.
  \end{aligned}
\end{equation}

In the following, assume that $2\varepsilon^{-\frac{1}{1-\delta}}\in \mathds N$ for simplicity. Otherwise, we can replace $2\varepsilon^{-\frac{1}{1-\delta}}$ with $\lfloor 2\varepsilon^{-\frac{1}{1-\delta}}\rfloor$ throughout the proof of Proposition \ref{uppertail-0}. 
we divide the interval $[-n,n]$ into $2\varepsilon^{-\frac{1}{1-\delta}}$ small intervals, denoted by $I_i$, $i\in [-\varepsilon^{-\frac{1}{1-\delta}},\varepsilon^{-\frac{1}{1-\delta}}-1]$, each with length $\varepsilon^{\frac{1}{1-\delta}}n$.

Assume that the event $B_{n,\varepsilon}$ occurs.
For odd $i\in [0,\varepsilon^{-\frac{1}{1-\delta}}-1]$, we define the internal resistance in the interval $I_i$ by
\begin{equation*}
R_i=R_{I_{i-1}\cup I_i\cup I_{i+1}}(I_{i-1},I_{i+1}).
\end{equation*}
From the independence of edges, we can see that $\{R_i\}_{\text{odd } i\in [0,\varepsilon^{-\frac{1}{1-\delta}}-1]}$ are independent (see Figure \ref{event-B} for an illustration).
In addition, by \cite[Theorem 1.1]{DFH25+}, there exists a constant $c_3=c_3(\beta)>0$ (depending only on $\beta$) such that
\begin{equation}\label{prob-R}
  \mathds{P}\left[R_i\geq c_3(\varepsilon^{\frac{1}{1-\delta}}n)^\delta\mid B_{n,\varepsilon}\right]\geq \frac{1}{2}.
\end{equation}

Next, we estimate the number of indices $i\in [0,\varepsilon^{-\frac{1}{1-\delta}}-1]$ for which the event defined in  \eqref{prob-R} occurs.
For that, define
\begin{equation*}
\eta_i=
  \begin{cases}
    1, & \mbox{if the event in \eqref{prob-R} occurs}, \\
    0, & \mbox{otherwise}.
  \end{cases}
\end{equation*}

Combining \eqref{prob-R}, the independence of $\{R_i\}_{\text{odd } i\in [0,\varepsilon^{-\frac{1}{1-\delta}}-1]}$ and the Chernoff bound, we have the following estimate.
\begin{lemma}\label{numb-good}
  For all $\beta>0$, we have
  \begin{equation*}
    \mathds{P}\left[\sum_{\text{odd } i\in [0,\varepsilon^{-\frac{1}{1-\delta}}-1]}\eta_i\leq \frac{1}{8}\varepsilon^{-\frac{1}{1-\delta}}\ | B_{n,\varepsilon} \right]\leq \exp\left\{-\frac{1}{32}\varepsilon^{-\frac{1}{1-\delta}}\right\}.
  \end{equation*}
  For simplicity, we will refer to the event $\{\sum_{\text{odd } i\in [0,\varepsilon^{-\frac{1}{1-\delta}}-1]}\eta_i\leq \frac{1}{8}\varepsilon^{-\frac{1}{1-\delta}}\}^c$ as $F_\varepsilon$.
\end{lemma}

We now present the

\begin{proof}[Proof of Proposition \ref{uppertail-0}]
Fix $n\in \mathds{N}$ and $\varepsilon\in (0,1)$.
It is easy to see that on the event $B_{n,\varepsilon}\cap F_\varepsilon$, for any unit flow $f$ in the \eqref{def-Rhat} (the red curves in Figure \ref{event-B}), it must pass through each small intervals $I_i$, $i\in [0,\varepsilon^{-\frac{1}{1-\delta}}-1]$ (see Figure \ref{event-B} for an illustration).
\begin{figure}[htbp]
\centering
\includegraphics[scale=0.6]{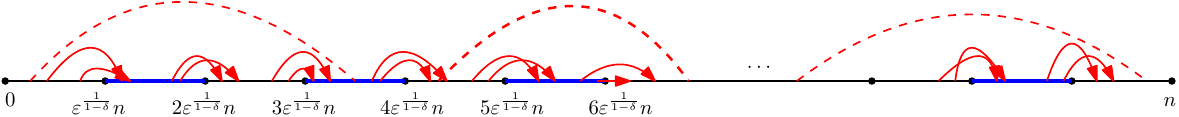}
\caption{The illustration for the event $B_{n,\varepsilon}$ and the independence of $\{R_i\}_{\text{odd } i\in [0,\varepsilon^{-\frac{1}{1-\delta}}-1]}$.}
\label{event-B}
\end{figure}
Consequently, we get that on the event $B_{n,\varepsilon}\cap F_\varepsilon$,
\begin{equation}\label{R-lb}
\widehat{R}([-n,0],[-n,n]^c)\geq \sum_{\text{odd } i\in [0,\varepsilon^{-\frac{1}{1-\delta}}-1]} R_i
\geq \frac{1}{8}\varepsilon^{-\frac{1}{1-\delta}}\cdot c_3(\varepsilon^{\frac{1}{1-\delta}}n)^\delta\geq \frac{c_3}{8} \varepsilon^{-1}n^\delta,
\end{equation}
where $c_3>0$ (depending only on $\beta$) is the constant defined in \eqref{prob-R}.
Similarly,  on the event $B_{n,\varepsilon}\cap F_\varepsilon$, we have
\begin{equation*}\label{R-lb2}
\widehat{R}([0,n],[-n,n]^c)\geq \frac{c_3}{8} \varepsilon^{-1}n^\delta.
\end{equation*}
Combining this with \eqref{R-lb}, \eqref{R-Rhat} and the fact $\max\{|\theta_1|,|\theta_2|\}\geq 1/2$, one has that on the event $B_{n,\varepsilon}\cap F_\varepsilon$,
\begin{equation*}
  R(0,[-n,n]^c)\geq \frac{c_3}{32} \varepsilon^{-1}n^\delta.
\end{equation*}
Therefore, by \eqref{prob-B} and Lemma \ref{numb-good}, we obtain that
\begin{equation*}
\begin{aligned}
  \mathds{P}\left[R(0,[-n,n]^c)\geq \frac{c_3}{32} \varepsilon^{-1}n^\delta\right]&\geq \mathds{P}\left[B_{n,\varepsilon}\cap F_\varepsilon\right]\\
  &=\mathds{P}\left[B_{n,\varepsilon}\right]\mathds{P}\left[ F_\varepsilon\ | B_{n,\varepsilon}\right]\geq \frac{1}{2}\exp\left\{-6\beta \varepsilon^{-\frac{1}{1-\delta}}\right\}.\quad \quad  \qedhere
  \end{aligned}
\end{equation*}
\end{proof}

It is worth mentioning here that 
the bound of the upper tail in Propostion~\ref{uppertail-0} is sharp according to the following 
upper bound. 
The proof is given 
in Section~\ref{better-up-tail} below.
\begin{proposition}\label{uppertail-up}
  For all $\beta>0$ and all $\eta\in(0,\frac{1}{1-\delta})$, there exists a constant $c_\eta>0$ {\rm(}depending only on $\beta$ and $\eta${\rm)} such that for all $n\in\mathds{N}$ and $\varepsilon\in(0,1)$ with $\varepsilon^{1/(1-\delta)}n\geq 1$, 
  \begin{equation*}
    \mathds{P}\left[\max_{x\in[-n,n]}R(0,x)\geq \varepsilon^{-1}n^\delta\right]\leq c_\eta^{-1}\exp\left\{-c_\eta\varepsilon^{-\eta}\right\}. 
  \end{equation*} 
 \end{proposition}

\medskip
At the end of the subsection, we provide the following estimate of the degree of a interval, which is also useful for the verification of Assumption \ref{assumpC}(1). 

\begin{lemma}\label{est-deg}
For all $\beta>0$, there exist constants $c_4>0$ and $M_0<\infty$ {\rm(}depending only on $\beta${\rm)} such that for all $n\in \mathds{N}$, 
  \begin{equation*}\label{degtail}
    \mathds{P}\left[\sum_{u\in[-n,n]}\mathrm{deg}(u)\geq M_0 n\right]\leq \exp\{-c_4 n\}.
  \end{equation*}
\end{lemma}
\begin{proof}
According to the independence of edges and the fact that $u\sim u\pm 1$ for all $u\in \mathds{Z}$, it can be checked that
\begin{equation}\label{stoc}
  \sum_{u\in[-n,n]}\mathrm{deg}(u)\leq_d 2\sum_{u\in [-n,n]}\sum_{v\in\mathds{Z}\setminus\{u,u\pm1\}}\xi_{uv}+4n,
\end{equation}
where $X\leq_d Y$ means that $X$ is stochastically dominated by $Y$, and  
$\{\xi_{uv}\}_{u\in[-n,n],\ v\in\mathds{Z}\setminus\{u,u\pm1\}}$ is a sequence of independent random variables with
\begin{equation*}
\mathds{P}\left[\xi_{uv}=1\right]=\mathds{P}\left[u\sim v\right]\quad \text{and}\quad \mathds{P}\left[\xi_{uv}=0\right]=1-\mathds{P}\left[u\sim v\right].
\end{equation*}
Thus, we can see that
\begin{equation*}
  \mathds{E}\left[\sum_{u\in [-n,n]}\sum_{v\in\mathds{Z}\setminus\{u,u\pm1\}}\xi_{uv}\right]= \widetilde{C}_1n,
\end{equation*}
where $\widetilde{C}_1<\infty$ is a constant depending only on $\beta$.
Combining this with the Chernoff bound and \eqref{stoc}, we obtain 
\begin{equation*}
  \mathds{P}\left[\sum_{u\in[-n,n]}\mathrm{deg}(u)\geq 4(\widetilde{C}_1+1)n\right]
  \leq \mathds{P}\left[\sum_{u\in [-n,n]}\sum_{v\in\mathds{Z}\setminus\{u,u\pm1\}}\xi_{uv}\geq 2\widetilde{C}_1n\right]\leq \e^{-\widetilde{C}_1n/3}.
\end{equation*}
Hence, the proof is complete.
\end{proof}

\subsection{Proof of Proposition~\ref{uppertail-up}}\label{better-up-tail}
This section is devoted to proving Proposition~\ref{uppertail-up}. 
To the end, for all $\beta>0$ and $n\in \mathds{N}$, we define  
\begin{equation*}
  \Lambda(n,\beta)=\max_{i,j\in [0,n)}\mathds{E}\left[R_{[0,n)}(i,j)\right]. 
\end{equation*}
For simplicity, we will denote $\Lambda(n,\beta)$ by $\Lambda(n)$ when there is no risk of confusion. Note that, from \cite[(7.4)]{DFH25+} we can see that $\Lambda(n)\asymp n^\delta$, where $\delta \in (0,1)$ is the exponent in Theorem \ref{eff-tightness}. 

Using $\Lambda(n)$, we first establish a suitable upper bound for moments of effective resistances. 
It is worth noting that estimates for the moments of resistances in the $\beta$-LRP model have already been provided in \cite[Lemma~3.3]{FH26}. 
The following lemma serves as an enhanced version of that. 
Since the proof shares similarities with that of \cite[Lemma~3.3]{FH26}, we will omit certain details to focus on the differences. 

\begin{lemma}\label{moment-est}
  For all $\beta>0$, there exists a constant $C_1,\,C_2<\infty$ 
  depending only on $\beta$ such that for all $r>0,\ n \in \mathds{N}$ and all $i,j \in [0,n)$,
  \begin{equation*}\label{moment-goal}
    \mathds{E}[R_{[0,n)}(i,j)^r]\leq C_1r^{C_2r}\Lambda(n)^r.
  \end{equation*}
\end{lemma}

\begin{proof}
We start with the case where $r\in \mathds{N}$. For any $n\in \mathds{N}$, define  
  \begin{equation*}
    \Gamma(n,r)=\max_{i,j\in[0,n)}\mathds{E}\left[R_{[0,n)}(i,j)^r\right].
  \end{equation*}
  For a fixed $m\in \mathds{N}$, from \cite[(3.2)-(3.5)]{FH26} we get the following iterative inequality: 
  \begin{equation}\label{induction}
    \Gamma(m^{n+1},2^k)\leq 2^{2^k}\widetilde{f}(\beta,m)\Gamma(m^n,2^k)+16m^2(8m)^{2^k} \left(\Gamma(m^n,2^{k-1})\right)^2
  \end{equation}
  for all $n,k\in\mathds{N}$, where 
  \begin{equation*}
    \widetilde{f}(\beta,m):=
      \begin{cases}
        2+\frac{20}{1-\beta}m^{1-\beta},\quad &\beta<1,\\
        12+8\log m,\quad &1\leq \beta\leq 2,\\
        22, \quad &\beta>2.
      \end{cases}
  \end{equation*}

  In the following, we aim to show that there exists a constant $M<\infty$ (depending only on $\beta$ and $m$) such that 
  \begin{equation}\label{hypothesis}
    \Gamma(m^n,2^k)\leq M^{k2^k}\Lambda(m^n)^{2^k}\quad \text{for all }n,k\in \mathds{N}. 
  \end{equation}
  For that, we will employ induction. Specifically, 
  By \cite[Lemma~3.3]{FH26}, we can see that \eqref{hypothesis} holds for $k=1$ and all $n\in\mathds{N}$. 
  Now assume that \eqref{hypothesis} holds for all $k'<k$ and all $n\in \mathds{N}$. 
  Applying this to \eqref{induction}, we get that 
  \begin{equation*}\label{induction-1}
    \Gamma(m^{n+1},2^k)\leq 2^{2^k}\widetilde{f}(\beta,m)\Gamma(m^n,2^k)+16m^2(8m)^{2^k} M^{(k-1)2^k}\Lambda(m^n)^{2^k}.
  \end{equation*}
  Furthermore, using a similar iteration as \cite[(3.6)-(3.7)]{FH26}, one has 
  \begin{equation*}
  \begin{aligned}
    \Gamma(m^{n+1},2^k)&\leq 16m^2(8m)^{2^k} M^{(k-1)2^k}\sum_{l=1}^n \left(2^{2^k}\widetilde{f}(\beta,m)\right)^{n+1-l}\Lambda(m^l)^{2^k}\\
    &\leq (\widetilde{C}_1m^2)^{2^k}M^{(k-1)2^k}\Lambda(m^{n+1})^{2^k}
    \end{aligned}
  \end{equation*}
  for some constant $\widetilde{C}_1<\infty$ depending only on $\beta$. 
  Thus, if we set $M=\widetilde{C}_1m^2$,  we can conclude that \eqref{hypothesis} holds for all $n\in\mathds{N}$, thereby completing the induction.

We now 
show that \eqref{hypothesis} holds for $\Gamma(N,r)$ for all $N\in \mathds{N}$ and $r\in \mathds{R}_+$. 
In fact, for general $N\in \mathds{N}$, let $p$ be the number such that $m^p\leq N<m^{p+1}$. 
Then for any $u,v\in [0,N)$ with $u<v$, we can find a sequence $\{i_k\}_{k=1}^K\subset \{1,2,\cdots,\lceil Nm^{-p}\rceil\}$ such that $|i_k-i_{k+1}|=1$, and $u\leq i_1m^p<\cdots<i_Km^p\leq v$. 
From this and the triangle inequality of the resistance, we get that for each $k\in \mathds{N}$, 
\begin{equation*}
\begin{aligned}
  \mathds{E}\left[R_{[0,N)}(u,v)^{2^k}\right]&\leq m^{2^k}\bigg[\mathds{E}\left[R_{[(i_1-1)m^p,i_1m^p)}(u,i_1m^p-1)^{2^k}+1\right]\\
  &\quad +\sum_{k=1}^K\mathds{E}\left[R_{[i_km^p,(i_k+1)m^p)}(i_km^p, (i_k+1)m^p-1)^{2^k}\right] +1\\
  &\quad + \mathds{E}\left[R_{[i_Km^p,(i_K+1)m^p)}(i_Km^p, v)^{2^k}\right]\bigg]. 
  \end{aligned}
\end{equation*}
Combining this with \eqref{hypothesis}, we conclude that \eqref{hypothesis} holds for $\Gamma(N,r)$ for all $N\in \mathds{N}$ and $r=2^1,2^2,\cdots$. 

Finally, for general $r\in \mathds{R}_+$, let $k\in \mathds{N}$ be the number such that $2^{k-1}<r\leq 2^k$. Then by the H\"older inequality and \eqref{hypothesis}, we have 
  \begin{equation*}
    \Gamma(N,r)\leq \Gamma(N,2^k)^{r/2^k}\leq M^{kr}\Lambda(N)^r\leq r^{\widetilde{C}_2r}\Lambda(N)^r
  \end{equation*}
  for some $\widetilde{C}_2<\infty$ depending only on $\beta$. Hence, the proof is complete. 
\end{proof}

With Lemma~\ref{moment-est} at hand,  we can get the following exponential moment of the effective resistances. 

\begin{lemma}\label{weak-mgf}
  For all $\beta>0$, there exists a constant $\eta>0$ depending only on $\beta$ such that
  \begin{equation}\label{eq-weak-mgf}
    \sup_{n\in\mathds{N}}\max_{i,j\in[0,n)}\mathds{E}\left[\exp\left\{\left(\frac{R_{[0,n)}(i,j)}{n^\delta}\right)^\eta\right\}\right]<\infty.
  \end{equation}
\end{lemma}
\begin{proof}
  For any $\eta>0$, by Lemma~\ref{moment-est} we get that 
  \begin{equation*}
      \mathds{E}\left[\exp\left\{\left(\frac{R_{[0,n)(i,j)}}{\Lambda(n)}\right)^\eta\right\}\right] = \mathds{E}\left[\sum_{r\geq 0}\frac{1}{r!}\left(\frac{R_{[0,n)(i,j)}}{\Lambda(n)}\right)^{\eta r}\right]\leq \sum_{r\geq 0} \frac{C_1(\eta r)^{C_2\eta r}}{r!}, 
  \end{equation*}
  where $C_1,\,C_2<\infty$ is the constant in Lemma \ref{moment-est} depending only on $\beta$.  By Stirling's formula $r!\approx \sqrt{2\pi n}(n/\e)^n$, 
   it is clear that the above series converge for some sufficiently small $\eta<C_1^{-1}$. Combining this with the fact $\Lambda(n)\asymp n^\delta$ (see \cite[(7.4)]{DFH25+}) gives the desired result. 
\end{proof}

From this lemma, we can derive the following upper tail estimate for the resistances. It is important to note that \cite[Lemma~6.4]{Baumler23a} provides a similar upper tail for the graph distance in the $\beta$-LRP model, our proof is primarily inspired by the approach used in its proof. 

\begin{lemma}\label{twopoint-diam}
  For all $\beta>0$, if \eqref{eq-weak-mgf} holds for some $\eta>0$, then there exists a constant $C_{\eta}>0$ such that for each $\varepsilon\in (0,1)$ and each $n\in\mathds{N}$,  
  \begin{equation}\label{uppertail-2}
    \mathds{P}\left[\max_{n'\in[0,n]}\max_{i,j\in[0,n']} R_{[0,n')}(i,j)>\varepsilon^{-1}n^\delta\right]\leq C_{\eta}\exp\{-\varepsilon^{-\eta}\}.
  \end{equation}
This implies that 
  \begin{equation}\label{weak-mgf-diam}
    \sup_{n\in\mathds{N}}\mathds{E}\left[\exp\left\{\left(\frac{\max_{i,j\in[0,n)} R_{[0,n)}(i,j)}{n^\delta}\right)^{\eta'}\right\}\right]<\infty\quad \text{for all }\eta'\in(0,\eta).
  \end{equation}
\end{lemma}

\begin{proof}
We start with the case $n=2^m$ for some $m\in \mathds{N}$. 
For $\varepsilon\in (0,1)$ and $l\in [0,m-1]$, let $F_l$ be the event that 
\begin{equation*}
  \max_{y\in [2^l,2^m)\cap 2^l\mathds{Z}}R_{[y-2^l,y)}(y-2^l,y)\leq c_\delta 2^{-\delta}\left(\frac{3}{4}\right)^{\delta(m-l)}\varepsilon^{-1}2^{\delta m}, 
\end{equation*}
where $c_\delta>0$ is a sufficiently small constant depending only on $\delta$. For each integer $x\in[0,2^m)$, let
\[
x_l=2^l\lfloor x/2^l\rfloor,\qquad 0\leq l\leq m.
\]
Then $x_0=x$, $x_m=0$, and $x_l-x_{l+1}\in\{0,2^l\}$. Applying the resistance triangle inequality along these sequences for the two endpoints, and using
$\sum_{l\geq0}(3/4)^{l\delta}<\infty$, shows that on $\bigcap_{l=0}^{m-1}F_l$ the maximum in \eqref{uppertail-2} is at most $\varepsilon^{-1}n^\delta$. Therefore, by a union bound, we get that 
\begin{equation}\label{diam-step-1}
\begin{aligned}
&\mathds{P}\left[\max_{n'\in[0,n]}\max_{i,j\in[0,n']} R_{[0,n')}(i,j)> \varepsilon^{-1}n^\delta\right] \leq \mathds{P}\left[\bigcup_{l=0}^{m-1}F_l^c\right]\\
&\leq \sum_{l=0}^{m-1}2^{m-l+1}\mathds{P}\left[R_{[0,2^l)}(0,2^l)\geq c_\delta 2^{-\delta}\left(\frac{3}{4}\right)^{\delta(m-l)}\varepsilon^{-1}2^{\delta m}\right].
\end{aligned}
\end{equation}

In addition, it follows from  Lemma \ref{weak-mgf} and Markov's inequality that for each $l\in [0,m-1]$, 
\begin{equation*}
\begin{aligned}
&\mathds{P}\left[R_{[0,2^l)}(0,2^l)> c_\delta 2^{-\delta}\left(\frac{3}{4}\right)^{\delta(m-l)}\varepsilon^{-1}2^{\delta m}\right]\\
&=\mathds{P}\left[\frac{R_{[0,2^l)}(0,2^l)}{2^{l\delta}}>c_\delta 2^{-\delta}\left(\frac{3}{2}\right)^{\delta(m-l)}\varepsilon^{-1}\right]\\
&\leq \mathds{E}\left[\exp\left\{\left(\frac{R_{[0,2^l)}(0,2^l)}{2^{l\delta}}\right)^\eta\right\}\right]
\exp\left\{-(c_\delta 2^{-\delta})^\eta\varepsilon^{-\eta}\left(\frac{3}{2}\right)^{\delta(m-l)\eta}\right\}\\
&\leq \widetilde{C}_1 \exp\left\{-(c_\delta 2^{-\delta})^\eta\varepsilon^{-\eta}\left(\frac{3}{2}\right)^{\delta(m-l)\eta}\right\}
\end{aligned}
\end{equation*}
for some $\widetilde{C}_1<\infty$ depending only on $\beta$ and $\eta$. Consequently, we obtain \eqref{uppertail-2} by applying this to \eqref{diam-step-1}. 

For general $n\in \mathds{N}$, let $m$ be the number such that $2^{m-1}\leq n<2^{m}$. It is clear that 
\begin{equation*}
  \begin{aligned}
  \mathds{P}\left[\max_{n'\in[0,n]}\max_{i,j\in[0,n']} R_{[0,n')}(i,j)> \varepsilon^{-1}n^\delta\right]
  &\leq \mathds{P}\left[\max_{n'\in[0,2^m]}\max_{i,j\in[0,n']} R_{[0,n')}(i,j)> \varepsilon^{-1}n^\delta\right]\\
  &\leq  \mathds{P}\left[\max_{n'\in[0,2^m]}\max_{i,j\in[0,n']} R_{[0,n')}(i,j)> (2^\delta\varepsilon)^{-1}2^{\delta m}\right]. 
  \end{aligned}
\end{equation*}
Hence, we can complete the proof by using \eqref{uppertail-2} with $n=2^m$.  
\end{proof}

Inspired by the argument used in the proof of \cite[Lemma~6.7]{Baumler23a} regarding graph distance, we can also enhance the value of $\eta'$ in \eqref{weak-mgf-diam} as follows. 

\begin{lemma}\label{improve-eta}
For $\beta>0$, assume that \eqref{eq-weak-mgf} holds for some $\eta>0$. Then we have 
  \begin{equation*}
    \sup_{n\in\mathds{N}}\mathds{E}\left[\exp\left\{\left(\frac{\max_{i,j\in[0,n)} R_{[0,n)}(i,j)}{n^\delta}\right)^{\eta'}\right\}\right]<\infty\quad \text{for all }\eta'\in(0,1+\delta\eta). 
  \end{equation*}
\end{lemma}

\begin{proof}
For each $n\in\mathds{N}$, define $a_0(n)=0$ and 
  \begin{equation*}
    a_{k+1}(n)=a_k(n)+2+\inf\left\{z\in\mathds{N}:R_{[a_k(n)+2,a_{k}(n)+2+z)}(a_k(n)+2,a_k(n)+2+z)>n^\delta\right\}
  \end{equation*}
  for all $k\geq 1$. 
  Note that $(a_{k+1}(n)-a_k(n))_{k\geq 0}$ are {\it i.i.d.}\ random variables since $(a_{k+1}(n)-a_k(n))_{k\geq 1}$ depends on edges in distinct intervals. 
  Moreover, for $\alpha\in(0,1)$, from the translation invariance of the $\beta$-LRP model and Lemma~\ref{twopoint-diam}, we have
  \begin{equation}\label{ak-tail}
    \begin{aligned}
      \mathds{P}[a_{1}(n)-a_0(n)\leq \alpha n]&=\mathds{P}\left[\exists z\leq \alpha n\text{ such that }R_{[0,z)}(0,z)>n^\delta\right]\\
      &\leq \mathds{P}\left[\max_{z\in[0,\alpha n]}\max_{i,j\in[0,z]}R_{[0,z)}(i,j)>\alpha^{-\delta}(\alpha n)^\delta\right]
      \leq C_{\eta}\exp\left\{-\alpha^{-\eta\delta}\right\}, 
    \end{aligned}
  \end{equation}
  where $C_\eta<\infty$ is the constant in Lemma \ref{twopoint-diam}.

  Now let $K_n=\inf\{k\geq 1:a_k(n)\geq n\}$ for all $n\in \mathds{N}$. 
  Then it follows that  
  \begin{equation*}
   R_{[0,n)}(0,n)\leq K_n n^\delta+2K_n\leq 3K_n n^\delta. 
  \end{equation*}
   Therefore, for any large integer $L\geq 1$, if $K_n>2L$,  there must exist at least $L$ integers $i\in\{1,2,\cdots,2L\}$ such that $a_{i}(n)-a_{i-1}(n)\leq L^{-1}n$. Otherwise, we can deduce that 
   \begin{equation*}
      n>a_{2L}=\sum_{i=1}^{2L}(a_i-a_{i-1})\geq L\cdot L^{-1}n=n. 
   \end{equation*}
   This leads to a contradiction.  
   As a result,
  \begin{equation*}
    \begin{aligned}
      \mathds{P}\left[R_{[0,n)}(0,n)>6Ln^\delta\right]&\leq \mathds{P}\left[K_n>2L\right]\\
      &\leq 2^{2L}\mathds{P}\left[a_i(n)-a_{i-1}(n)\leq L^{-1}n\right]^L\\
      &\leq 2^{2L}C_{\eta}^L\exp\left\{-L^{1+\eta\delta}\right\}, 
    \end{aligned}
  \end{equation*}
  where the last inequality uses \eqref{ak-tail}. This implies that 
  \begin{equation*}
    \sup_{n\in\mathds{N}}\mathds{E}\left[\exp\left\{\left(\frac{R_{[0,n)}(0,n)}{n^\delta}\right)^{\widetilde{\eta}}\right\}\right]<\infty
  \end{equation*}
  for all $\widetilde{\eta}\in (0,1+\delta \eta)$. 
  Furthermore, by the monotonicity of the internal resistance and the translation invariance of the $\beta$-LRP model, one has 
  \begin{equation}\label{supmax}
  \begin{aligned}
  \sup_{n\in\mathds{N}}\max_{i,j\in[0,n)}\mathds{E}\left[\exp\left\{\left(\frac{R_{[0,n)}(i,j)}{n^\delta}\right)^{\widetilde{\eta}}\right\}\right]&\leq \sup_{n\in\mathds{N}}\max_{i,j\in[0,n):i<j}\mathds{E}\left[\exp\left\{\left(\frac{R_{[i,j)}(i,j)}{|j-i|^\delta}\right)^{\widetilde{\eta}}\right\}\right]\\
  &=\sup_{n\in\mathds{N}}\mathds{E}\left[\exp\left\{\left(\frac{R_{[0,n)}(0,n)}{n^\delta}\right)^{\widetilde{\eta}}\right\}\right]<\infty.
  \end{aligned}
  \end{equation}
  That is, \eqref{eq-weak-mgf} holds for $\widetilde{\eta}\in(0,1+\delta\eta)$. 
  Finally, we can get the desired result by using the similar argument in the proof of Lemma \ref{twopoint-diam}. 
\end{proof}

In the final part of this subsection, we present the

\begin{proof}[Proof of Proposition~\ref{uppertail-up}]
  Let $f(\eta)=1+\delta\eta$. From Lemma~\ref{weak-mgf}, we can see that \eqref{eq-weak-mgf} holds for some small $\eta>0$. 
  Furthermore, using \eqref{supmax} gives that \eqref{eq-weak-mgf} holds for each $\eta'\in(0,f(\eta))$. 
  Iterating such argument we get that \eqref{eq-weak-mgf} holds for all $\eta\in(0,f^{k}(\eta))$ and  $k\in\mathds{N}$. 
  Note that for any fixed $\eta\in(0,1/(1-\delta))$, we have $f^{k}(\eta)\to 1/(1-\delta)$ as $k\to\infty$. Thus for each $\eta\in(0,1/(1-\delta))$, we have \eqref{eq-weak-mgf}. Combining this with Lemma \ref{twopoint-diam} yields the desired result. 
\end{proof}

\subsection{Proof of Propositions \ref{Vn-lu}(1) and \ref{prop-assumpC2}}\label{subsect-largeball}

This subsection is devoted to proving Propositions~\ref{prop-assumpC2} and \ref{Vn-lu}(1). 
We start by constructing, with exponentially small but controlled probability, a resistance network across a prescribed interval. 
The idea is to partition the interval into blocks of size $m$ and connect these blocks
hierarchically by long edges, while controlling the resistance within the relevant subintervals using Proposition~\ref{uppertail-up}.

More precisely, for each $m, s\in \mathds N$, let $K=2^s$ and divide the interval $[-Km,Km)$ into small intervals with length $m$, denoted by 
$$
I_{m,j}=[(-K+j)m,\, (-K+j+1)m)\quad\text{for }0\leq j<2K. 
$$ 
  By Proposition~\ref{uppertail-up}, we can see that there exists a constant $D_*\geq 1$, depending only on $\beta$, such that for any finite interval $I$,
\begin{equation}\label{Dstar}
\mathds P\left[\operatorname{diam}(I;R)\leq D_*|I|^\delta\right]\geq 1/2.
\end{equation}
The following lemma gives the resulting resistance bound.

\begin{lemma}\label{lem-sepcont}
For all $\beta>0$, there exist constants  $c_5\geq 1$ and $C_3<\infty$, depending only on $\beta$, such that for every $s\in \mathds N$ and every $m\geq m_0(K):=\lceil2\left(1+\log_2(2K)\right)^{2/\delta}\rceil$, 
\begin{equation*}
  \mathds P\left[\max_{x\in [-Km,Km)}R_{[-Km,Km)}(0,x)\leq c_5m^\delta\right]\geq {\rm e}^{-C_3K}, 
\end{equation*}
where $K:=2^s$. 
\end{lemma}

\begin{figure}[htbp]
\centering
\includegraphics[scale=0.46]{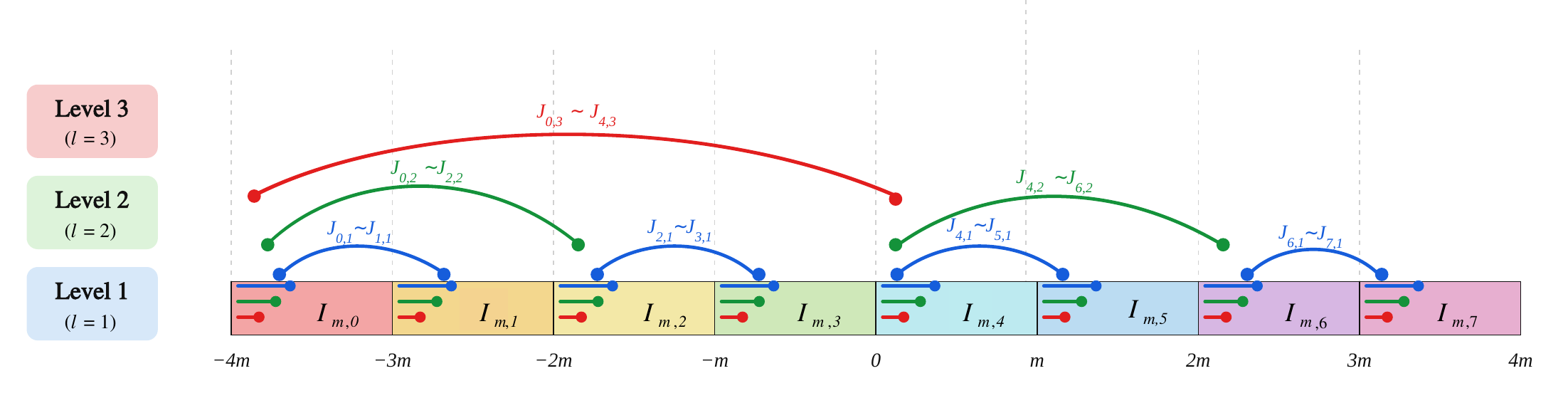}
\caption{The illustration of the hierarchical tree construction for $s=2$. The blocks $I_{m,0},\ldots, I_{m,7}$ form the leaves of the construction. Within each block, the blue, green and red segments represent $J_{j,1},\ J_{j,2}$ and $J_{j,3}$, respectively, while arcs of the corresponding colors indicate the level-1, level-2 and level-3 tree edges. The port intervals are nested as $J_{j,3}\subseteq J_{j,2}\subseteq J_{j,1}$.}
\label{N8}
\end{figure}

\begin{proof}
For fixed $s\in \mathds N$ and $m\geq m_0(K)$, let 
\begin{equation*}
  J_{j,l}=[(-K+j)m,\, (-K+j)m+t_l),\quad t_l=\left\lfloor\frac{m}{(l+1)^{2/\delta}}\right\rfloor\quad \text{for all }1\leq l\leq s+1. 
\end{equation*}
By the definition of $m_0(K)$ and the fact that $l\leq s+1<1+\log_2(2K)$,
we have
\begin{equation}\label{eq-tl}
\frac{m}{2(l+1)^{2/\delta}}
\leq t_l
\leq
\frac{m}{(l+1)^{2/\delta}}
\quad\text{and}\quad 
m^\delta\geq s+1.
\end{equation}

For $1\leq l\leq s+1$ and $0\leq q<2K/2^l$, let $A_{l,q}$ denote the event that there is an edge between $J_{q2^l,l}$ and $J_{q2^l+2^{l-1},l}$.
Moreover, let $B_{l,q}$ be the event 
$$
\operatorname{diam}\left(J_{q2^l,l};R\right)\vee
\operatorname{diam}\left(J_{q2^l+2^{l-1},l};R\right)
\leq D_*t_l^\delta,
$$  
and let $F_j$ be the event that 
$$
\operatorname{diam}(I_{m,j};R)\leq D_*m^\delta\quad \text{for }0\leq j<2K,  
$$ 
where $D_*<\infty$ is the constant defined in \eqref{Dstar}. 
Denote  
\begin{equation}\label{def-FKm}
G_{K,m}
=
\left(
\bigcap_{l=1}^{s+1}\,
\bigcap_{q=0}^{2K/2^l-1} \left(A_{l,q}\cap B_{l,q}\right)
\right)
\cap
\left(
\bigcap_{j=0}^{2K-1} F_j. 
\right). 
\end{equation}
Notice that all events in the right hand-side of \eqref{def-FKm} are increasing with respect to edges. Therefore, according to \eqref{Dstar}, 
FKG inequality and the translation invariance of $\beta$-LRP model, we get that 
\begin{equation}\label{probFKm}
\mathds P\left[G_{K,m}\right]
\geq \prod_{l=1}^{s+1}\prod_{q=0}^{2K/2^l-1}\mathds P\left[A_{l,q}\right]\mathds P\left[B_{l,q}\right]
\prod_{j=0}^{2K-1}\mathds P\left[F_j\right]\geq 2^{-(6K-2)}\prod_{l=1}^{s+1}\left(\mathds P\left[A_{l,0}\right]\right)^{2K/2^l}.
\end{equation}
Additionally, it follows from \eqref{prob-pij} that for every $1\leq l\leq s+1$, 
\[
\mathds P\left[A_{l,0}\right]
=1-\mathds P\left[A_{l,0}^c\right]
=1-\exp\left\{-\beta\int_{J_{0,l}}\int_{J_{2^{l-1},l}}\frac{1}{|u-v|^2}\d u\d v\right\}\geq \widetilde c_14^{-l}(l+1)^{-4/\delta}
\]
for some $\widetilde c_1>0$ depending only on $\beta$. Applying this to \eqref{probFKm}, we get that
\begin{equation*}
-\log\mathds P[G_{K,m}]
\leq (6K-2)\log 2+
\sum_{l=1}^{s+1}\frac{2K}{2^l}
\left(\log\frac{1}{\widetilde c_1}+l\log 4+\frac4\delta\log(l+1)\right)
\leq \widetilde C_1K
\end{equation*}
for some $\widetilde C_1<\infty$ depending only on $\beta$, where we may take $\widetilde c_1\leq1$. Consequently,
\begin{equation}\label{eq-probFKmfinal}
\mathds P[G_{K,m}]\geq\e^{-\widetilde C_1K}.
\end{equation}

It remains to show that
\[
G_{K,m}\subseteq \left\{\max_{x\in [-Km,Km)}R_{[-Km,Km)}(0,x)\leq c_1m^\delta\right\}.  
\]
Indeed, on $G_{K,m}$, choose one present edge for each pair of ports specified by
the events $A_{l,q}$.
These edges form the hierarchical tree on the block
labels $\{0,\ldots,2K-1\}$ with root block $I_{m,0}$. In addition, at level $l$, the tree connects
$J_{q2^l,l}$ to $J_{q2^l+2^{l-1},l}$,  
$0\leq q<2K/2^l$. 
See Figure~\ref{N8} for the case $s=2$.

Fix $x\in I_{m,j}$ for some $j$. The levels of the tree edges along the path from $I_{m,0}$ to $I_{m,j}$ are strictly decreasing. Since $t_l$ is nonincreasing in $l$, the endpoints of the incoming and outgoing edges in an intermediate block lie in the port corresponding to the outgoing edge. The resistance between them is therefore bounded by $D_*t_l^\delta$ at that level. The same bound applies to the initial motion from $-Km$, while the final motion to $x$ costs at most $D_*m^\delta$. Hence, by the resistance triangle inequality, Rayleigh monotonicity and \eqref{eq-tl},
\begin{align*}
R_{[-Km,Km)}(-Km,x)
&\leq D_*m^\delta+D_*\sum_{l=1}^{s+1}t_l^\delta+(s+1)\nonumber\\
&\leq D_*m^\delta\left(1+\sum_{l=1}^{\infty}(l+1)^{-2}\right)+m^\delta
\leq (2D_*+1)m^\delta.\label{eq-Ra0x}
\end{align*}
Applying this also to $x=0$ and using the resistance triangle inequality again, we obtain
\begin{equation}\label{eq-GKm}
G_{K,m}\subseteq
\left\{\max_{x\in[-Km,Km)}R_{[-Km,Km)}(0,x)\leq(4D_*+2)m^\delta\right\}.
\end{equation}
Together with \eqref{eq-probFKmfinal}, this proves the desired estimate with $c_5:=4D_*+2$ and $C_3:=\widetilde C_1$.
\end{proof}

The preceding construction controls the resistance from the origin
throughout the core interval $[-Km,Km)$. 
For the multiscale argument,
we also need to isolate this core from the exterior by a region of
large effective resistance. 
We therefore enlarge the core by adding
buffer intervals of width $Lm$ on both sides. 
The next lemma shows
that, with a comparable exponential probability, the resistance inside
the enlarged interval remains controlled while the resistance from the
core to its exterior is bounded from below.

\begin{lemma}\label{lem-singlesc}
For all $\beta>0$, there exist constants $c_6\geq 8c_5$, $C_4<\infty$ and integers
$L,K_0\geq1$, depending only on $\beta$, such that the following holds.
For $s\in \mathds N$ with $K=2^s\geq K_0$ and every $m\geq m_0(K)$, 
with probability at least $\mathrm e^{-C_4K}$, 
\begin{equation*}\label{eq-event12}
\begin{aligned}
&\max_{x\in[-Km,Km)}
R_{[-(K+L)m,(K+L)m)}(0,x)
\leq c_5m^\delta,\\
&\max_{x\in [-(K+L)m,(K+L)m)}
R_{[-(K+L)m,(K+L)m)}(0,x)
\leq c_6m^\delta,
\end{aligned}
\end{equation*}
and
\begin{equation*}\label{eq-event3}
R([-Km,Km),\, [-(K+L)m,(K+L)m)^c)
\geq20c_5m^\delta.
\end{equation*}
\end{lemma}

\begin{proof}
Fix $L\geq1$. Let 
\[
J_-=[-(K+L)m,-Km)
\quad\text{and}\quad 
J_+=[Km,(K+L)m)
\]
be the left and right buffer intervals, respectively, and let
\[
S_-=[-Km,-Km+2Lm)
\quad\text{and}\quad 
S_+=[Km-2Lm,Km)
\]
be the corresponding boundary portions of the core interval
$[-Km,Km)$.

Assume that $K\geq4L+4$, and let $F_{K,m}$ be the event that none
of the following edges is present (see Figure \ref{noedgep} below):
\begin{equation}\label{noedge}
\begin{aligned}
&[-Km,Km)\sim [-(K+L)m,(K+L)m)^c,\quad
J_-\sim[-Km,Km)\setminus S_-,\\
&J_+\sim[-Km,Km)\setminus S_+,
\quad
J_-\sim J_+.
\end{aligned}
\end{equation}
Here $I\sim J$ means that $I$ and $J$ are connected by edges directly. 
From the independence of edges, we can see that $F_{K,m}$ is independent of
$G_{K,m}$ defined in \eqref{def-FKm}.
For simplicity, we denote by $A_i$ the $i$-th event in \eqref{noedge} for $i=1,\cdots,4$. 
Now by \eqref{prob-pij} and the independence of edges again, we have 
\begin{equation}\label{prob-FKm}
\mathds P\left[F_{K,m}\right]
=\prod_{i=1}^4\mathds{P}\left[A_i^c\right]
\geq \widetilde c_1(K+1)^{-2\beta}
\end{equation}
for some $\widetilde c_1>0$ depending only on $\beta$. 
\begin{figure}[htbp]
\centering
\includegraphics[scale=0.6]{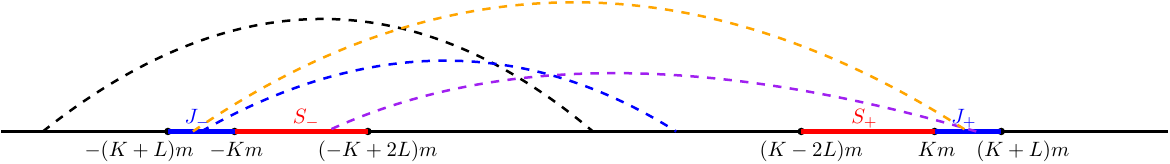}
\caption{The illustration of \eqref{noedge}. The black dashed curve represent the absence of edges between $[-Km,Km)$ and $[-(K+L)m,(K+L)m)^c$, while the blue (resp.\ purple) dashed curve represent the absence of edges between $J_-$ (resp.\ $J_+$) and $[-Km,Km)\setminus S_-$ (resp.\ $[-Km,Km)\setminus S_+$). Additionally, the orange dashed curve represent the absence of edges between $J_-$ and $J_+$. Therefore, for the flow from $[-Km,Km)$ to $[-(K+L)m,(K+L)m)^c$, it must pass through the interval $J_-\cup J_+$.}
\label{noedgep}
\end{figure}

On $F_{K,m}$, we wire $[-Km,Km)$ to a single source vertex and
$[-(K+L)m,(K+L)m)^c$ to a single sink vertex. 
By \eqref{noedge}, the resulting network consists of two parallel channels through $J_-$ and $J_+$. Denote their effective
resistances by $\mathcal R_-$ and $\mathcal R_+$, respectively. 
Then we can see that 
\begin{equation}\label{parallel}
R([-Km,Km),\,[-(K+L)m,(K+L)m)^c)
=
\left(\mathcal R_-^{-1}+\mathcal R_+^{-1}\right)^{-1}.
\end{equation}

We next bound the resistance of the two channels. Let us start with $\mathcal R_+$.
Conditional on $F_{K,m}\cap G_{K,m}$, the edges internal to $J_+$, from $J_+$ to $S_+$, and from $J_+$ to
$[-(K+L)m,(K+L)m)^c$ retain their original independent laws. 
Therefore, it follows from Theorem \ref{eff-tightness} that there exists a constant $\widetilde c_2>0$ (depending only on $\beta$) such that 
\begin{equation}\label{eq-mathcalR+}
\mathds P\left[
\mathcal R_+\geq \widetilde c_2(Lm)^\delta
\,\middle|\,
F_{K,m}\cap G_{K,m}
\right]
\geq\frac{15}{16}.
\end{equation}
By the symmetry, we can see that the same estimate holds for $\mathcal R_-$. 

In addition, by the independence of edges again and Proposition~\ref{uppertail-up}, 
there is a constant $\widetilde C_1<\infty$, depending only on $\beta$, such that
\[
\mathds P\left[
\operatorname{diam}(J_\pm;R)\leq \widetilde C_1(Lm)^\delta
\,\middle|\,
F_{K,m}\cap G_{K,m}
\right]
\geq\frac{15}{16}.
\]
Combining this with \eqref{eq-GKm}, \eqref{parallel} and \eqref{eq-mathcalR+}, we obtain that conditionally on 
$F_{K,m}\cap G_{K,m}$, with probability at least $3/4$,
\begin{equation*}\label{eq-KLm}
R([-Km,Km),\, [-(K+L)m,(K+L)m)^c)
\geq\frac{\widetilde c_2}{2}(L m)^\delta
\end{equation*}
and 
\begin{equation}\label{eq-RD}
\max_{x\in [-(K+L)m,(K+L)m)}
R_{[-(K+L)m,(K+L)m)}(0,x)
\leq
\left(c_5+1+\widetilde C_1L^\delta\right)m^\delta.
\end{equation}
In particular, on $G_{K,m}$, Lemma~\ref{lem-sepcont} and Rayleigh
monotonicity imply 
\begin{equation}\label{eq-DKm}
\max_{x\in[-Km,Km)}
R_{[-(K+L)m,(K+L)m)}(0,x)
\leq c_5m^\delta.
\end{equation}

Finally, choose $L$ sufficiently large such that $\widetilde c_5L^\delta\geq40c_1$,  
and set
\[
c_6=
\max\left\{8c_5,\;c_5+1+\widetilde C_1L^\delta\right\}
\quad\text{and} \quad 
K_0=4L+4.
\]
Combining Lemma~\ref{lem-sepcont}, \eqref{prob-FKm}, \eqref{eq-RD} and \eqref{eq-DKm}, we obtain the probability of the event in the lemma is bounded below by 
\[
\frac34\,\widetilde c_1(K+1)^{-2\beta}\mathrm e^{-C_3K}
\geq
\mathrm \e^{-\widetilde C_2K},
\]
for all $K\geq K_0$ and some $\widetilde C_2<\infty$ depending only on $\beta$. 
\end{proof}

For $H\geq1$, let
\begin{equation}\label{lrp:K-choice}
K(H)=
2^{\left\lceil
\log_2\!\left(\max\left\{K_0,\,(4c_5)^{d_{f}}H\right\}\right)
\right\rceil}.
\end{equation}
Then it is easy to check that $K(H)\asymp_\beta  H$.
Consequently, Lemma~\ref{lem-singlesc}, applied with $K=K(H)$,
holds with probability at least
\begin{equation}\label{eq-probC3}
\exp\left\{-C_5H\right\}
\end{equation}
for some $C_5<\infty$ depending only on $\beta$.

We now let $(H_k)_{k\geq1}$ be any deterministic sequence with $H_k\geq1$,
and set $K_k=K(H_k)$ for simplicity. 
Fix $\eta>0$. Starting from $T_0:=0$, define recursively
\begin{equation}\label{def-mkTk}
\begin{split}
m_k
&=
\max\left\{
m_0(K_k),
\left\lceil
4(k+1)^{\eta/2}(T_{k-1}+1)
\right\rceil
\right\},\\
T_k
&=
\left\lceil
2(k+1)^{1+\eta/2}(K_k+L)m_k
\right\rceil.
\end{split}
\end{equation}
By the above definitions, we can see that 
\begin{equation*}\label{eq-seper}
[-T_{k-1},T_{k-1}]
\subset
[-K_km_k,K_km_k)
\subset
[-(K_k+L)m_k,(K_k+L)m_k)
\subset
[-T_k,T_k],
\end{equation*}
and
\begin{equation}\label{eq-seper2}
\frac{T_{k-1}+1}{m_k}
\leq
\frac{1}{4(k+1)^{\eta/2}}.
\end{equation}

For each $k\geq1$, contract $[-T_{k-1},T_{k-1}]$ to a single vertex $*$ in $G[[-T_k,T_k]]$, retaining parallel edges and discarding self-loops, where $G[I]$ represents the subgraph of $\beta$-LRP model within $I$. 
For convenience, we use the same notation for the remaining vertices after contraction.
For $I=[-(K_k+L)m_k,(K_k+L)m_k)$ or $[-T_k,T_k]$, write $\widetilde R_I$ for effective resistance in
the contraction graph.

Let $\mathcal A_k$ be the event that
\begin{equation}\label{def-mathcalAk}
\begin{split}
&\max_{x\in[-K_km_k,K_km_k)}
\widetilde R_{[-(K_k+L)m_k,(K_k+L)m_k)}(*,\,x)
\leq c_5m_k^\delta,\\
&\max_{x\in[-(K_k+L)m_k,(K_k+L)m_k)}
\widetilde R_{[-(K_k+L)m_k,(K_k+L)m_k)}(*,\,x)
\leq c_6m_k^\delta,\\
&\widetilde R_{[-T_k,T_k]}
\left(
[-K_km_k,K_km_k),\,
[-T_k,T_k]\setminus
[-(K_k+L)m_k,(K_k+L)m_k)
\right)
\geq20c_5m_k^\delta.
\end{split}
\end{equation}
where $c_5,c_6>0$ are constants defined in Lemmas \ref{lem-sepcont} and \ref{lem-singlesc}, respectively.

The contraction above separates the randomness used at different
scales. This makes the events
$(\mathcal A_k)$ independent. The next lemma combines this observation
with the single-scale estimate and then transfers the resulting
resistance bounds from the contracted network back to the original
network. 
In particular, it shows that a successful trial at scale $k$ produces the desired event $\mathcal H$ uniformly over the admissible values of $\lambda$.

\begin{lemma}\label{lem-probAk}
The events $(\mathcal A_k)_{k\geq1}$ are independent and satisfy
\begin{equation*}
\mathds P\left[\mathcal A_k\right]\geq \exp\left\{-C_5H_k\right\}, 
\end{equation*}
where $C_5<\infty$ is the constant defined as in \eqref{eq-probC3}. 
Moreover, almost surely, for all sufficiently large $k$, on
$\mathcal A_k$,
\begin{equation}\label{eq-unifksuc}
\mathcal H(16c_1\lambda m_k^\delta,\,\lambda,\,H_k)
\quad\text{occurs for every}\quad
\lambda\geq\frac{c_6}{4c_5}\geq2.
\end{equation}
\end{lemma}

\begin{proof}
From the definition of $\mathcal A_k$ in \eqref{def-mathcalAk}, we see that the event $\mathcal A_k$ depends only on edges $\langle u,v\rangle$ satisfying
\[
u,v\in[-T_k,T_k]
\quad\text{and}\quad 
\{u,v\}\not\subset[-T_{k-1},T_{k-1}],
\quad
|u-v|>1.
\]
Note that for different $k$, these sets of edges are disjoint. 
Indeed, if $j<k$,
then every edge determining $\mathcal A_j$ has both endpoints in
$
[-T_j,T_j]\subset[-T_{k-1},T_{k-1}],
$ 
whereas such edges are excluded from those determining $\mathcal A_k$. 
Hence
$(\mathcal A_k)_{k\geq1}$ are independent by the independence of edges in the model. 

To bound the probability of $\mathcal A_k$, we add fresh independent
edges with probability \eqref{prob-pij} inside $[-T_{k-1},T_{k-1}]$, and extend the resulting
network to $\mathds Z$. 
By using Lemma~\ref{lem-singlesc} and \eqref{eq-probC3} with $K=K_k$ and $m=m_k$, we claim that  
at least
\begin{equation}\label{eq-probAk}
\mathds P\left[\mathcal A_k\right]\geq \exp\left\{-C_5H_k\right\}.
\end{equation}
Indeed restricting the network to $[-T_k,T_k]$ can only increase the set
resistance in the third inequality of \eqref{def-mathcalAk},
while the first two inequalities involve only
\[
[-(K_k+L)m_k,(K_k+L)m_k)
\subset[-T_k,T_k].
\]
In addition, contracting $[-T_{k-1},T_{k-1}]$ can only decrease the point
resistances, and it leaves the set resistance unchanged since
\[
[-T_{k-1},T_{k-1}]
\subset[-K_km_k,K_km_k),
\]
which is already wired when that resistance is computed. Thus the
event in Lemma~\ref{lem-singlesc} implies $\mathcal A_k$. This proves \eqref{eq-probAk}. 

We next transfer these estimates back to the original network.
From Proposition~\ref{uppertail-up} with $\eta=1$ there and \eqref{eq-seper2} we have
\begin{equation}\label{eq-interr}
\begin{split}
\mathds P\left[A_k\right]&:=
\mathds P\left[
\operatorname{diam}\left([-T_{k-1},T_{k-1}];\,R\right)>c_5m_k^\delta
\right]\\
&\leq
\widetilde C\exp\left\{-\widetilde{C}^{-1}
\left(
\frac{2T_{k-1}+1}{m_k}
\right)^{-1}\right\}
\leq \widetilde C\exp\left\{-\widetilde{C}^{-1}
(k+1)^{\eta/2}\right\}, 
\end{split}
\end{equation}
where $0<\widetilde C<\infty$ is a constant depending only on $\beta$. 
Moreover, using $p_{i,j}\leq \widetilde C_\beta|i-j|^{-2}$ by \eqref{prob-pij} and $T_k\geq2(K_k+L)m_k$, 
we obtain
\begin{equation*}\label{lrp:external-error}
\begin{aligned}
\mathds P\left[B_k\right]&:=\mathds P\left[
[-(K_k+L)m_k,(K_k+L)m_k)
\sim [-T_k,T_k]^c
\right]\\
&\leq
\widetilde C_\beta\frac{(K_k+L)m_k}{T_k}
\leq
\widetilde C'_\beta(k+1)^{-1-\eta/2}, 
\end{aligned}
\end{equation*}
where constants $\widetilde C_\beta, \widetilde C'_\beta<\infty$ depend only on $\beta$. 
Combining this with \eqref{eq-interr} by taking $p>2/(\eta\delta)$, we get that 
\begin{equation*}
\sum_{k=1}^\infty\mathds P\left[A_k\cup B_k\right]
\leq\widetilde \sum_{k=1}^\infty \widetilde C\exp\left\{-\widetilde{C}^{-1}(k+1)^{\eta/2}\right\}
+\widetilde C'_\beta\sum_{k=1}^\infty(k+1)^{-1-\eta/2}<\infty. 
\end{equation*}
Hence, by the first Borel--Cantelli lemma,
almost surely, for all sufficiently large $k$, $A^c_k\cap B^c_k$ occurs. That is, 
\begin{equation}\label{eq-ioABc}
\operatorname{diam}\left([-T_{k-1},T_{k-1}];\,R\right)
\leq c_2m_k^\delta \quad 
\text{and}\quad 
[-(K_k+L)m_k,(K_k+L)m_k)
\not\sim[-T_k,T_k]^c.
\end{equation}

On $\mathcal A_k$, whenever \eqref{eq-ioABc} holds, then according to the first two inequalities in the definition of $\mathcal A_k$ in \eqref{def-mathcalAk} and the triangle inequality of resistance, we have 
\begin{equation}\label{eq-Kkmk}
\begin{aligned}
&\max_{x\in[-K_km_k,K_km_k)}
R_{[-(K_k+L)m_k,(K_k+L)m_k)}(0,x)\\
&\leq \max_{x\in[-K_km_k,K_km_k)}
\widetilde R_{[-(K_k+L)m_k,(K_k+L)m_k)}(*,\,x)
+\operatorname{diam}\left([-T_{k-1},T_{k-1}];\,R\right)
\leq2c_5m_k^\delta,
\end{aligned}
\end{equation}
and
\begin{equation}\label{eq-Kkmk2}
\begin{aligned}
&\max_{x\in[-(K_k+L)m_k,(K_k+L)m_k)}
R_{[-(K_k+L)m_k,(K_k+L)m_k)}(0,x)\\
&\leq \max_{x\in[-(K_k+L)m_k,(K_k+L)m_k)}
\widetilde R_{[-(K_k+L)m_k,(K_k+L)m_k)}(*,\,x)
+\operatorname{diam}\left([-T_{k-1},T_{k-1}];\,R\right)\\
&\leq
(c_6+c_5)m_k^\delta
\leq2c_6m_k^\delta.
\end{aligned}
\end{equation}
In addition, if 
$
[-(K_k+L)m_k,(K_k+L)m_k)
\not\sim[-T_k,T_k]^c,
$ 
since
$
[-T_{k-1},T_{k-1}]
\subset[-K_km_k,K_km_k),
$
we get 
\begin{equation}\label{eq-Rlower}
\begin{aligned}
&R\left(
[-K_km_k,K_km_k),
[-(K_k+L)m_k,(K_k+L)m_k)^c
\right)\\
&=\widetilde R_{[-T_k,T_k]}
\left(
[-K_km_k,K_km_k),\,
[-T_k,T_k]\setminus
[-(K_k+L)m_k,(K_k+L)m_k)
\right)
\geq20c_5m_k^\delta.
\end{aligned}
\end{equation}

Now let $\lambda\geq c_6/(4c_5).$ 
By \eqref{eq-Kkmk} we have 
\begin{equation}\label{eq-Kkmksub}
[-K_km_k,K_km_k)
\subset B(0,4c_5m_k^\delta),
\end{equation}
while \eqref{eq-Kkmk2} implies 
\[
[-(K_k+L)m_k,(K_k+L)m_k)
\subset B(0,16c_5\lambda m_k^\delta).
\]
Hence, by the monotonicity of resistance and \eqref{eq-Rlower}, we arrive at 
\begin{align*}
R\left(
0,B(0,16c_5\lambda m_k^\delta)^c
\right)
&\geq
R\left(
[-K_km_k,K_km_k),
[-(K_k+L)m_k,(K_k+L)m_k)^c
\right)\\
&\geq20c_5m_k^\delta
>
\frac{16c_5\lambda m_k^\delta}{2\lambda}.
\end{align*}

Finally, from \eqref{eq-Kkmksub}, the definition of $K_k$ in \eqref{lrp:K-choice} and the fact that 
$\bm V(0,r)\geq 2|B(0,r)|$ in $\beta$-LRP model, we see that 
\begin{align*}
\bm V(0,4c_5m_k^\delta)
&\geq2\,|[-K_km_k,K_km_k)|
=4K_km_k
\geq4H_k(4c_5)^{d_{f}}m_k
=4H_k(4c_5m_k^\delta)^{d_{f}}
>H_k(4c_5m_k^\delta)^{d_{f}},
\end{align*}
where we used $\delta d_{f}=1$.
These estimates prove \eqref{eq-unifksuc} simultaneously for
every $\lambda\geq c_6/(4c_5)$.
\end{proof}

Additionally, to convert these above scales into the asymptotic growth of the radii appearing in Proposition~\ref{prop-assumpC2}, we next determine the growth rate of the recursively defined sequence $(m_k)$. 
The
following deterministic estimate will be used after selecting the successful indices.

\begin{lemma}\label{lrp:scale-growth}
Suppose that $\log K_k=o(\log k)$ and that, for all sufficiently large $k$,
\[
m_k=
\left\lceil
4(k+1)^{\eta/2}(T_{k-1}+1)
\right\rceil.
\]
Then we have $\log m_k=(1+\eta)k\log k\,(1+o(1)).$ 
In particular, these assumptions hold if
$ 
H_k=1+\varepsilon\log(k+{\rm e})
$ 
for any fixed $\varepsilon>0$.
\end{lemma}

\begin{proof}
For all sufficiently large $k$, it follows from \eqref{def-mkTk} that 
\[
\log m_k-\log m_{k-1}
=(1+\eta)\log k+\log(K_{k-1}+L)+O(1).
\]
Summing over $k$ and using $\log K_k=o(\log k)$, we obtain
\[
\log m_k=(1+\eta)k\log k+o(k\log k).
\]
If $H_k=1+\varepsilon\log(k+\mathrm e)$, then $K_k\asymp\log k$ and
\begin{equation*}\label{eq-growm0}
m_0(K_k)=O\left((1+\log\log k)^{2/\delta}\right).
\end{equation*}
On the other hand, \eqref{def-mkTk} implies $m_k>4(k+1)^{\eta/2}m_{k-1}$. Thus the second term in the maximum defining $m_k$ eventually exceeds $m_0(K_k)$, and the assumptions hold.
\end{proof}

With the above lemmas at hand, we now present the 

\begin{proof}[Proof of Proposition \ref{prop-assumpC2}] 
Throughout the proof, we take $\lambda_0=\frac{c_6}{4c_5},$ and $C_0:=C_5$, where $c_5,\,c_6>0$ and $C_5<\infty$ are constants defined as in Lemmas \ref{lem-sepcont}, \ref{lem-singlesc} and  \ref{lem-probAk}, respectively. 
Set 
\[
H_k=1+\frac{\theta}{C_0}\log(k+\e).
\]
and let $(m_k)$ be defined by \eqref{def-mkTk} with this choice of $(H_k)$. 

From Lemma~\ref{lem-probAk}, we see that the events
$(\mathcal A_k)_{k\geq1}$ are independent and
\[
\mathds P\left[\mathcal A_k\right]
\geq
\exp\left\{-C_5H_k\right\}
=
\e^{-C_0}(k+\e)^{-\theta}.
\]
Hence, there exists a constant $\widetilde c_1>0$ (depending only on $\beta$) such that for all sufficiently large $n$,
\begin{equation*}\label{eq-Akunif}
\mathds P\left[
\bigcap_{k=\lfloor n^\alpha\rfloor}^{\lfloor(n+1)^\alpha\rfloor-1}
\mathcal A_k^c
\right] 
\leq
\exp\left\{
-\sum_{k=\lfloor n^\alpha\rfloor}^{\lfloor(n+1)^\alpha\rfloor-1}
\mathds P\left[\mathcal A_k\right]
\right\}
\leq
\exp\left\{-\widetilde c_1 n^{\alpha(1-\theta)-1}\right\}.
\end{equation*}
Since $\alpha(1-\theta)-1>0$ by \eqref{def-alpha},  
\[
\sum_{n\geq1}
\mathds P\left[
\bigcap_{k=\lfloor n^\alpha\rfloor}^{\lfloor(n+1)^\alpha\rfloor-1}
\mathcal A_k^c
\right]
<\infty.
\]
Hence, by the first Borel--Cantelli lemma, almost surely only finitely
many of these intervals contain no index $k$ for which $\mathcal A_k$
occurs. Equivalently, almost surely every sufficiently large interval
$
\{\lfloor n^\alpha\rfloor,\ldots,
\lfloor(n+1)^\alpha\rfloor-1\}
$ 
contains an index $k$ for which $\mathcal A_k$ occurs.
Together with Lemma~\ref{lem-probAk}, this allows us to choose a finite random integer $J\geq0$ such that, for every $n\geq1$,
there exists an index $k_n$ satisfying
\begin{equation}\label{cond-kn}
\lfloor(n+J)^\alpha\rfloor
\leq k_n
<
\lfloor(n+J+1)^\alpha\rfloor
\end{equation}
and
$
\mathcal H
\left(
16c_5\lambda m_{k_n}^{\delta},
\lambda,
H_{k_n}
\right)
$ 
occurs for every $\lambda\geq\frac{c_6}{4c_5}.$
We now choose $k_n$ to be the smallest such index and set
\begin{equation}\label{eq-defradii}
r(n)=16c_5\lambda(n)m_{k_n}^{\delta}.
\end{equation}

Note that $(k_n)$ is strictly increasing, $(m_k)$ is strictly
increasing and $\lambda(n)$ is nondecreasing. Thus $(r(n))$ defined in \eqref{eq-defradii} is also strictly
increasing. 
Moreover, by \eqref{cond-kn} and the fact $g(n)\ge 1$, one has 
\begin{align*}
H_{k_n}&=1+\frac{\theta}{C_0}\log(k_n+\e)
\geq
1+\frac{\theta}{C_0}
\log(\lfloor n^\alpha\rfloor+\e)
=h(n),\\
\lambda(n)&=\lambda_0g(n)
\geq\lambda_0
=\frac{c_6}{4c_5}.
\end{align*}
Thus applying Lemma~\ref{lem-probAk} with 
$\lambda=\lambda(n)$, we have 
$\mathcal H(r(n),\lambda(n),H_{k_n})$ occurs for every $n\geq1$.
Furthermore, since $H_{k_n}\geq h(n)$, by decreasing the volume threshold from
$H_{k_n}$ to $h(n)$, we get that 
$
\mathcal H(r(n),\lambda(n),h(n))
$
occurs for every $n\geq1$.

It remains to verify the two asymptotic statements. By \eqref{cond-kn}, we have 
\[
\frac{k_n}{n^\alpha}\rightarrow1
\quad\text{and}\quad 
\frac{\log k_n}{\log n}\rightarrow\alpha\qquad \text{as }n\rightarrow \infty.
\]
In addition, combining this with  Lemma~\ref{lrp:scale-growth}, we get 
\begin{equation}\label{est-mk}
\log m_{k_n}
=
(1+\eta)k_n\log k_n
+
o(k_n\log k_n)
=
\alpha(1+\eta)n^\alpha\log n
+
o(n^\alpha\log n).
\end{equation}
Since \eqref{def-functg} implies
\[
\log\lambda(n)
=
\log\lambda_0+\log g(n)
=
o(\log n),
\]
from this, \eqref{eq-defradii} and \eqref{est-mk} we obtain
\[
\log r(n)
=
\alpha(1+\eta)\delta n^\alpha\log n
+
o(n^\alpha\log n).
\]
This proves \eqref{lrp:target-speed}.

Finally, by the definition of $h(n)$, we see that 
\[
h(n)
=
1+\frac{\theta}{C_0}
\log(\lfloor n^\alpha\rfloor+\e)
=
\frac{\alpha\theta}{C_0}\log n+O(1),
\]
while
$\lambda(n)^{1+d_{f}}
=
\lambda_0^{1+d_{f}}
g(n)^{1+d_{f}}.$
Consequently,
\[
\frac{h(n)}{\lambda(n)^{1+d_{f}}}
=
\frac{\alpha\theta}
{C_0\lambda_0^{1+d_{f}}}
\frac{\log n}{g(n)^{1+d_{f}}}
\,(1+o(1))
\rightarrow\infty
\]
by \eqref{def-functg}. 
This proves
\eqref{lrp:target-ratio}.
\end{proof}

In the remainder of this subsection, we prove Proposition \ref{Vn-lu}(1) by using Proposition \ref{prop-assumpC2} as follows. 

\begin{proof}[Proof of Proposition \ref{Vn-lu}(1)]
Fix $\eta,\theta,\alpha$ and $g$ satisfying \eqref{def-alpha} and \eqref{def-functg}, and let $(r(n))_{n\geq1}$ be the sequence in Proposition \ref{prop-assumpC2}. Define 
\[
b_n=\frac{r(n)}{4\lambda(n)}.
\]
Then from Proposition \ref{prop-assumpC2} and the definition of $\mathcal H(r,\lambda,h)$ in \eqref{def-mathcalH}, we get that almost surely, for every $n\geq1$,
\begin{equation}\label{eq-vtail}
\bm V(0,b_n)\geq h(n)b_n^{1/\delta}.
\end{equation}
Moreover, \eqref{def-functg} and \eqref{lrp:target-speed} imply
\[
\log b_n=\alpha(1+\eta)\delta n^\alpha\log n\, (1+o(1))
\quad\text{and}\quad 
\log\log b_n=\alpha\log n\, (1+o(1)). 
\]
Since $h(n)=(\alpha\theta/C_0)\log n(1+o(1))$, it follows \eqref{eq-vtail} that, for every $\alpha_1\in(0,1)$,
\begin{equation*}\label{eq:vol-limsup}
\limsup_{n\to\infty}
\frac{\bm V(0,b_n)}{b_n^{1/\delta}(\log\log b_n)^{\alpha_1}}
=\infty\qquad\mathds P\text{-a.s.}
\end{equation*}
Replacing $b_n$ by $\lceil b_n\rceil$ and using the monotonicity of $\bm V(0,\cdot)$, we complete the proof. 
\end{proof}

\subsection{Proof of Proposition \ref{Vn-lu}(2) and Theorem \ref{osc-LRP}}\label{subsect-proofThm} 

In this subsection, we primarily complete the proofs of Proposition \ref{Vn-lu}(2) and Theorem \ref{osc-LRP} by verifying that the $\beta$-LRP model satisfies the conditions in Assumptions \ref{assumpA}, \ref{assumpB} and \ref{assumpC}. 

We start with estimating the lower tail of resistance. 

\begin{lemma}\label{lem-rlow-hklu}
  For all $\beta>0$, there exist constants $c_7, q_1>0$ {\rm(}both depending only on $\beta${\rm)} such that  for all $x_0\in\mathds{Z}$, $r,\lambda\geq 1$ and all $\varepsilon\leq 1/(2\lambda)$.
  \begin{equation*}
    \mathds{P}\left[R(x,B(x_0,r)^c)\geq r/\lambda\quad \text{for all}\quad x\in B(x_0,\varepsilon r)\right]\geq 1-c_7\lambda^{-q_1}.
  \end{equation*}
\end{lemma} 
\begin{proof}
  From the translation invariance of the $\beta$-LRP model, assume that $x_0=0$ without loss of generality. 
  Similarly to \eqref{eq:resttri}, we 
  see that 
  if there exists $x\in B(0,\varepsilon r)$ such that $R(x,B(0,r)^c)< r/\lambda$, then 
  $$R(0,B(0,r)^c)\leq R(0,x)+R(x,B(0,r)^c)\leq \varepsilon r+r/\lambda <2r/\lambda.$$
  Thus it suffices to show that there exist constants $\widetilde{c}_1,\widetilde{q}>0$ depending only on $\beta$ such that
  \begin{equation}\label{r-low-goal}
    \mathds{P}\left[R(0,B(0,r)^c)<2r/\lambda\right]\leq \widetilde{c}_1 \lambda^{-\widetilde{q}}
  \end{equation}
  holds for each $r,\lambda\geq 1$. To this end, let $n=\lambda^{-1/(2\delta)}r^{1/\delta}$.  
 Note that if $R(0,B(0,r)^c)<2r/\lambda$ and $[-n,n]\subset B_r(0)$,  then we have  $R(0,[-n,n]^c)<2r/\lambda$. Therefore, 
  \begin{equation}\label{r-low-1}
    \mathds{P}\left[R(0,B(0,r)^c)<2r/\lambda\right]\leq \mathds{P}\left[R(0,[-n,n]^c)<2r/\lambda\right]+\mathds{P}\left[[-n,n]\nsubseteq B_r(0)\right].
  \end{equation}
  
  In the following, we aim to estimate the two probabilities in the right-hand side of \eqref{r-low-1}. 
  It follows from \cite[Corollary~1.2]{DFH25+} that there exist constants $q_0,\widetilde{c}_2>0$ such that for each $r,\lambda\geq 1$,
  \begin{equation}\label{r-low-2}
    \mathds{P}\left[R(0,[-n,n]^c)<2r/\lambda\right]=\mathds{P}\left[\frac{R(0,[-n,n]^c)}{n^\delta}<2\lambda^{-1/2}\right]\leq \widetilde{c}_2\lambda^{-q_0}.
  \end{equation}
  In addition, it is clear that on the event $[-n,n]\nsubseteq B_r(0)$, 
  $$
  \max_{x\in[-n,n]}R(0,x)\geq r.
  $$ 
  Consequently, from Proposition~\ref{uppertail-up} with $\eta=1$, there exist constants $\widetilde{C}, \widetilde{c}_3>0$ such that
  \begin{equation*}\label{r-low-3}
  \begin{aligned}
    \mathds{P}\left[[-n,n]\nsubseteq B_r(0)\right]&\leq\mathds{P}\left[\max_{x\in[-n,n]}R(0,x)\geq r\right] \\
    &=\mathds{P}\left[\max_{x\in[-n,n]}R(0,x)\geq \lambda^{1/2}n^\delta \right]
    \leq \widetilde{C}\exp\{-\widetilde{c}_3\lambda^{1/2}\}.
    \end{aligned}
  \end{equation*}
  Applying this and \eqref{r-low-2}  to \eqref{r-low-1} yields \eqref{r-low-goal}. 
\end{proof}

We also need the following lower tail for $\bm V(\cdot,\cdot)$ with the help of Proposition~\ref{uppertail-up}. 

\begin{lemma}\label{V-lower}
  For all $\beta>0$ and all  $\eta\in(0,\frac{1}{1-\delta})$, we have 
  \begin{equation*}
    \mathds{P}\left[\bm V(x,r)\leq \lambda^{-1}r^{1/\delta}\right]\leq 2c_\eta^{-1}\exp\left\{-c_\eta \lambda^{\delta\eta}\right\}
  \end{equation*}
   for all $x\in\mathds{Z}$ and all $r, \lambda\geq 1$,  where $c_\eta>0$ is the constant in Proposition \ref{uppertail-up} depending only on $\beta$ and $\eta$. 
\end{lemma}
\begin{proof}
  From \cite[(3.14)]{FH26}, we have
  \begin{equation*}
    \mathds{P}\left[\bm V(x,r)\leq \lambda^{-1}r^{1/\delta}\right]\leq 2\mathds{P}\left[\max_{x\in[0,\lambda^{-1}r^{1/\delta})}R(0,x)\geq r\right].
  \end{equation*}
  Combining this with Proposition~\ref{uppertail-up}, we can complete the proof. 
\end{proof}

The following upper tail of the volume $\bm V(\cdot,\cdot)$ is from \cite{FH26}. 
\begin{lemma}[{\cite[Lemma 3.6]{FH26}}]\label{up-V}
For all $\beta>0$, there exist constants $c_8,q_2>0$ {\rm(}both depending only on $\beta${\rm)} such that for any $x\in \mathds{Z}$ and any $r,\lambda\geq 1$, 
\begin{equation*}
  \mathds{P}\left[\bm V(x,r)>\lambda r^{1/\delta}\right]\leq c_8\lambda^{-q_2}. 
\end{equation*}
\end{lemma}

We next verify Assumption \ref{assumpC}(1). This will also imply Proposition \ref{Vn-lu}(2).
To this end, we introduce several events related to resistance and its volume growth.

Recall that $d_f=1/\delta$, and the events $H_1,H_2,H_3$ are defined in Definition~\ref{def-I-1}. For simplicity, we will set $c_3=c_4=1$ in Definition~\ref{def-I-1}. By Lemma \ref{V-lower}, it is clear that for any $\eta\in (0,\frac{1}{1-\delta})$, 
\begin{equation}\label{prob-H3}
  \mathds{P}\left[H_3(r,\lambda, u_2)^c\right]\leq 2c_\eta^{-1}\exp\left\{-c_\eta\left(4^{1/\delta}u_2\right)^{\eta\delta}\right\},  
\end{equation}
where $c_\eta>0$ is the constant in Proposition \ref{uppertail-up} depending only on $\beta$ and $\eta$.

To estimate the probability of the event $H_1(r,\lambda)$, we introduce the following event concerning the lower tail of the resistance.
\begin{definition}\label{def-H1}
 For $r>0$ and $\lambda\geq 1$, let $H_1'(r,\lambda)$ be the event that $R(0,B(0,r)^c)>2r/\lambda$.  
\end{definition}
According to the first paragraph in the proof of Lemma \ref{lem-rlow-hklu}, we can see that  $H_1'(r,\lambda)\subset H_1(r,\lambda)$.  
Additionally, we decompose the event $H_2(r,u_1)$ into the following two events. 
\begin{definition}\label{def-H21}
For $r>0$ and $u_1\geq 1$, let $H_{2,1}(r,u_1)$ be the event that 
\begin{equation*}
  R\left(0,[-r^{d_f}/(M_0u_1),\,r^{d_f}/(M_0u_1)]^c\right)\geq r, 
\end{equation*} 
 and let $H_{2,2}(r,u_1)$ be the event that 
 \begin{equation*}
   \sum_{z\in [-r^{d_f}/(M_0u_1),\,r^{d_f}/(M_0u_1)]}\text{deg}(z)\leq r^{d_f}/u_1. 
 \end{equation*}
 Here $M_0<\infty$ is the constant in Lemma \ref{est-deg} depending only on $\beta$.
 \end{definition}
 By the definition of volume $\bm V(\cdot,\cdot)$, one has $H_{2,1}(r,u_1)\cap H_{2,2}(r,u_1)\subset H_{2}(r,u_1)$. Furthermore, by Lemma \ref{est-deg}, we get that 
 \begin{equation}\label{prob-H22}
   \mathds{P}\left[H_{2,2}(r,u_1)^c\right]\leq \exp\left\{-c_4r^{d_f}/u_1\right\}, 
 \end{equation}
 where $c_4>0$ is the constant in Lemma \ref{est-deg} depending only on $\beta$. Moreover, it follows from Proposition \ref{uppertail-0} that 
 \begin{equation*}\label{prob-H21}
   \mathds{P}\left[H_{2,1}(r,u_1)\right]\geq \exp\left\{-c_2\left(M_0u_1\right)^{\frac{\delta}{1-\delta}}\right\},  
 \end{equation*}
  where $c_2>0$ is the constant in Proposition \ref{uppertail-0} depending only on $\beta$.

Now we turn to the 
\begin{proof}[Verification of Assumption \ref{assumpC}(1)] 
Due to the translation invariance of the $\beta$-LRP model, we only need to consider the case where $x=0$. We start with some notations which will be used repeatedly in the proof. Recall that $\delta \in (0, 1)$ is the exponent of resistance as defined in Theorem \ref{eff-tightness}. 
Throughout the proof, we will select $\eta\in(0,\frac{1}{1-\delta})$ and sufficiently small $\varepsilon\in(0,1)$, both depending only on $\beta$, such that
    \begin{equation}\label{eta}
      \alpha:=\frac{3(1-\varepsilon)(1-\delta)}{\delta}-\frac{1}{\delta\eta}\left(2+\frac{1}{1+\delta}\right)>0.
    \end{equation}
 Next, we define $\overline{\eta}=(\eta+(1-\delta)^{-1})/2\in(\eta,(1-\delta)^{-1})$ and $\overline{\varepsilon}=1-\varepsilon/2$, and take $\iota\in (0,1)$ sufficiently small satisfying 
 \begin{equation}\label{def-iota}
 \iota<\frac{\alpha}{4+6\delta^{-1}}\wedge \frac{(1-\varepsilon)(1-\delta)}{\overline{\varepsilon}}. 
 \end{equation}
 Fix $M>1$ depending only on $\beta$. All estimates below involving these scales are for sufficiently large $n$.
 For $n\in\mathds{N}$, let 
    \begin{equation*}\label{r-n}
    r_n=M^{n^2}\quad\text{and}\quad \lambda_n=(\log (\e-1+n))^{\iota}. 
    \end{equation*}
    We also denote 
    \begin{equation}\label{v1-n}
      h_2(n)=\widetilde{c}_1(\log (\e+n))^{\frac{(1-\varepsilon)(1-\delta)}{\delta}}/M_0\quad \text{and}\quad 
      h_3(n)=\widetilde{c}_2(\log (\e+n))^{1/(\eta\delta)},   
    \end{equation}
    where $\widetilde{c}_1,\widetilde{c}_2>0$ are normalizing constants such that $h_2(1)=h_3(1)=1$. 
    Then $\lambda,h_2,h_3\in\mathbf{SV}$. Recalling Assumption \ref{assumpC}(1), let 
    \begin{equation*}\label{Gn}
      G(n)=\frac{h_2(n)^3}{\lambda_n^{3d_f+2}h_3(n)^{2+\frac{d_f}{1+d_f}}}\in\mathbf{SV}.
    \end{equation*} 
    It is clear that $G(n)\asymp(\log(\e+n))^{\alpha_*}$. Here $\alpha_*:=\alpha-\iota(3/\delta+2)>0$, and we use \eqref{eta} and the fact that $d_f=1/\delta$ to simplify $G$. 
    With the choice  
    \begin{equation*}\label{Tn}
      T(n):=\left(\frac{r_n}{\lambda_n}\right)^{1+d_f}h_3(n)^{-1},  
    \end{equation*}
    we have $\Phi_2(n):=G(T^{-1}(n))\asymp(\log \log n)^{\alpha_*}$ as $n\to\infty$ and $\Phi_2\in {\bf SV}$. 
       
    Now it suffices to prove Assumption \ref{assumpC}(1)  for such $r_n, \lambda, h_2, h_3$, that is, to show that  $\mathds{P}$-a.s.\  $\{H_1(r_n,\lambda_n)\cap H_2(r_n,h_2(n))\cap H_3(r_n,\lambda_n,h_3(n))\text{ i.o.}\}$ occurs. 
    Furthermore, by Definition~\ref{def-H21}, one has 
    \begin{equation*}
      H_{2,1}(r_n,h_2(n))\cap H_{2,2}(r_n,h_2(n))\subseteq H_2(r_n,h_2(n)). 
    \end{equation*}
    Combining this with the fact $H_1'(r_n,\lambda_n)\subseteq H_1(r_n,\lambda_n)$, we only need to prove that $\mathds{P}$-a.s.\   $\{H_1'(r_n,\lambda_n)\cap H_{2,1}(r_n,h_2(n))\cap H_{2,2}(r_n,h_2(n))\cap H_3(r_n,\lambda_n,h_3(n))\text{ i.o.}\}$ occurs.

    We first claim that there are finitely many of the events $(H_{2,2}(r_n,h_2(n))\cap H_3(r_n,\lambda_n,h_3(n)))^c$ occurring. Indeed, according to \eqref{prob-H3} with $\eta=\overline{\eta}$, \eqref{prob-H22}, \eqref{v1-n} and the fact $\overline{\eta}=(\eta+1/(1-\delta))/2>\eta$, one has 
      \begin{align*}
        &\sum_{n\geq 1} \mathds{P}\left[H_{2,2}(r_n,h_2(n))^c\cup H_3(r_n,\lambda_n,h_3(n))^c\right]\\
        &\leq \sum_{n\geq 1} \left(\exp\left\{-c_4r_n^{d_f}h_2(n)^{-1}\right\}+
        2c_{\overline{\eta}}^{-1}\exp\left\{-c_{\overline{\eta}}\left(4^{d_f}h_3(n)\right)^{\delta\overline{\eta}}\right\}\right)\\
        &\leq \sum_{n\geq1}\left(\exp\left\{-c_4r_n^{d_f}h_2(n)^{-1}\right\}+2c_{\overline{\eta}}^{-1}\exp\left\{-\widetilde{c}_3(\log n)^{\overline{\eta}/\eta}\right\}\right)<\infty
      \end{align*}
      for some $\widetilde{c}_3>0$, where $c_4,c_{\overline{\eta}}>0$ are the constants (both depending only on $\beta$) defined in Lemma \ref{est-deg} (see also \eqref{prob-H22}) and Proposition \ref{uppertail-up}, respectively.  
    Combining this with the Borel-Cantelli lemma yields the claim. Therefore, we only need to prove that $\mathds{P}$-a.s.\ $\{H_1'(r_n,\lambda_n)\cap H_{2,1}(r_n,h_2(n))\text{ i.o.}\}$ occurs in the following. 
    
To achieve this, we need some new notations. Recall that $\overline{\varepsilon}=1-\varepsilon/2$ and $d_f=1/\delta$. For any $n\in \mathds{N}$, define   
\begin{equation}\label{def-an3}
\begin{aligned}
a_n&=h_2(n)^{-1}r_n^{d_f}/M_0 =\widetilde{c}_1^{-1}(\log (\e+n))^{-\frac{(1-\varepsilon)(1-\delta)}{\delta}}r_n^{1/\delta},\\ 
a_{n-1/2}&=\widetilde{c}_1^{-1}(\log (\e+n))^{-\frac{(1-\varepsilon)(1-\delta)}{\delta}}M^{(n-1/2)^2/\delta}
\end{aligned}
\end{equation} 
and 
\begin{equation*}
  w_n=(\log(\e+ n))^{-\frac{(1-\varepsilon)(1-\delta)}{\delta\overline{\varepsilon}}}r_n^{1/\delta}\in (a_{n-1/2},a_{n}).
\end{equation*}
Let $\mathcal{F}_n$ be the $\sigma$-field generated by edges with at least one endpoint contained in $[-a_n, a_n]$,  and for convenience, denote $\mathcal{F}_0=\emptyset$.  
    
In the following, we introduce two new events that are decompositions of the event $H_1'(r_n,\lambda_n)$. 
For any $n\in \mathds{N}$, let $H_{1,1}(r_n)$ be the event that 
\begin{equation*}\label{def-H11}
  \max_{z\in [-w_n,w_n]}R_{[-w_n,w_n]}(0,z)\leq r_n
\end{equation*}
and let   $H_{1,2}(r_n,\lambda_n)$ be the event that 
\begin{equation*}\label{def-H12}
  R(0,[-w_n,w_n]^c)> \frac{2r_n}{\lambda_n}. 
\end{equation*}   
By the definition of $H_1'$ (see Definition~\ref{def-H1}), it follows that 
\begin{equation}\label{dcp-H1}
H_{1,1}(r_n)\cap H_{1,2}(r_n,\lambda_n)\subseteq H_{1}'(r_n,\lambda_n).
\end{equation}
It is worth emphasizing that the events $H_{1,1}(r_n)$ and $H_{1,2}(r_n,\lambda_n)$ depends only on edges with at least one endpoint contained in $[-w_n,w_n]\subset [-a_n,a_n]$. This implies that $H_{1,1}(r_n), H_{1,2}(r_n,\lambda_n)\in \mathcal{F}_n$ for all $n\in\mathds{N}$. 
Moreover, by the definition of $H_{2,1}$ in Definition \ref{def-H21}, we also have $H_{2,1}(r_n,h_2(n))\in \mathcal{F}_n$ for all $n\in\mathds{N}$. 
Additionally, for any $n\in \mathds{N}$, let $H_n$ be the event that 
$[-a_{n-1},a_{n-1}]\nsim [-a_{n-1/2},a_{n-1/2}]^c$, and let $\widetilde{H}_n$ be the event that 
$$\max_{z\in [-a_{n-1},a_{n-1}]}R_{[-a_{n-1},a_{n-1}]}(0,z)\leq r_n/\lambda_n.$$ 
It is clear that $H_n,\widetilde H_n\in\mathcal F_{n-1}$. By Lemma \ref{twopoint-diam} with $\eta=1$ and \eqref{def-an3}, there exist constants $\widetilde C_1,\widetilde C_2<\infty$ depending only on $\beta$ such that
\begin{equation*}
\begin{aligned}
\sum_{n\geq1}\mathds P\left[H_n^c\cup\widetilde H_n^c\right]
&\leq\widetilde C_1\sum_{n\geq1}\frac{a_{n-1}}{a_{n-1/2}}
+\widetilde C_2\sum_{n\geq1}
\exp\left\{-\widetilde C_2^{-1}
\left(\frac{r_n}{\lambda_na_{n-1}^{\delta}}\right)\right\}
<\infty.
\end{aligned}
\end{equation*}
Here the finitely many initial terms may be bounded separately. By the Borel--Cantelli lemma, $H_n\cap\widetilde H_n$ occurs for all sufficiently large $n$ almost surely. Moreover, 
\begin{equation}\label{prob-HH}
\begin{aligned}
   &\I_{H_n\cap \widetilde{H}_n}\mathds{P}\left[H_{1,1}(r_n)\cap H_{1,2}(r_n,\lambda_n)\cap H_{2,1}(r_n,h_2(n))\mid \mathcal{F}_{n-1}\right]\\
  &\geq  \I_{H_n\cap \widetilde{H}_n}\left(\mathds{P}\left[H_{1,2}(r_n,\lambda_n)\cap H_{2,1}(r_n,h_2(n))\mid \mathcal{F}_{n-1}\right]- \mathds{P}\left[H_{1,1}(r_n)^c \mid \mathcal{F}_{n-1}\right]\right). 
  \end{aligned}
\end{equation}

Next, we will estimate the probabilities in the last inequality of \eqref{prob-HH} one by one. 
Note that on the event $H_n\cap \widetilde{H}_n$, by the monotonicity of the resistance, we have 
  \begin{equation*}
   H_{1,2}(r_n,\lambda_n)\supseteq \left\{ R([-a_{n-1/2},a_{n-1/2}], [-w_n,w_n]^c)\geq \frac{2r_n}{\lambda_n}\right\}=:  \widetilde{H}_{1,2}(r_n,\lambda_n)
 \end{equation*}
and 
\begin{equation*}
H_{2,1}(r_n,h_2(n))\supseteq \left\{ R([-a_{n-1/2},a_{n-1/2}], [-a_n,a_n]^c)\geq r_n\right\}=: \widetilde{H}_{2,1}(r_n,h_2(n)). 
\end{equation*}
Moreover, on the event $H_n$, it is easy to check that the resistances 
\begin{equation*}
R([-a_{n-1/2},a_{n-1/2}], [-w_n,w_n]^c)\quad \text{and}\quad R([-a_{n-1/2},a_{n-1/2}],[-a_{n},a_{n}]^c)
\end{equation*}  
both are independent of $\mathcal{F}_{n-1}$.  
Therefore, by the triangle inequality of the resistance and the FKG inequality for the two decreasing events below, we get that 
\begin{equation}\label{prob-H12H21}
  \begin{aligned}
  &\I_{H_n\cap \widetilde{H}_n}\mathds{P}\left[ H_{1,2}(r_n,\lambda_n) \cap H_{2,1}(r_n,h_2(n)) \mid \mathcal{F}_{n-1}\right]\\
  &\geq\I_{H_n\cap \widetilde{H}_n}\mathds{P}\left[\widetilde{H}_{1,2}(r_n,\lambda_n)\cap \widetilde{H}_{2,1}(r_n,h_2(n))\right]  \\
  &\geq \I_{H_n\cap \widetilde{H}_n}\mathds{P}\left[R(0,[-w_n,w_n]^c)\geq \frac{3r_n}{\lambda_n},\ R(0,[-a_n,a_n]^c)\geq 2r_n,\ \max_{z\in [-a_{n-1/2},a_{n-1/2}]} R(0,z)\leq \frac{r_n}{\lambda_n}\right]\\
  &\geq \I_{H_n\cap \widetilde{H}_n}\left(\mathds{P}\left[R(0,[-w_n,w_n]^c)\geq \frac{3r_n}{\lambda_n},\ R(0,[-a_n,a_n]^c)\geq 2r_n\right]\right.\\
  &\qquad \qquad \qquad \qquad \qquad\qquad \qquad\qquad \qquad \left. -\mathds{P}\left[\max_{z\in [-a_{n-1/2},a_{n-1/2}]} R(0,z)>\frac{r_n}{\lambda_n}\right]\right)\\
  &\geq \I_{H_n\cap \widetilde{H}_n}\left( \mathds{P}\left[R(0,[-w_n,w_n]^c)\geq \frac{3r_n}{\lambda_n}\right]\mathds{P}\left[R(0,[-a_n,a_n]^c)\geq 2r_n\right]\right.\\
  &\qquad \qquad \qquad \qquad \qquad\qquad \qquad\qquad \qquad\left.-\mathds{P}\left[\max_{z\in [-a_{n-1/2},a_{n-1/2}]} R(0,z)>\frac{r_n}{\lambda_n}\right]\right). 
  \end{aligned}
\end{equation}
Additionally, by the definitions of $a_n,w_n$ and $\lambda_n$, we have 
\begin{equation*}
\frac{r_n}{\lambda_nw_n^\delta}
\asymp(\log(\e+n))^{\frac{(1-\varepsilon)(1-\delta)}{\overline\varepsilon}-\iota}
\qquad\text{and}\qquad
\frac{r_n}{a_n^\delta}
\asymp(\log(\e+n))^{(1-\varepsilon)(1-\delta)}.
\end{equation*}
Thus by Proposition \ref{uppertail-0}, we get that 
\begin{equation*}\label{prob-H12H21-1}
\begin{aligned}
&\mathds P\left[R(0,[-w_n,w_n]^c)\geq\frac{3r_n}{\lambda_n}\right]
\mathds P\left[R(0,[-a_n,a_n]^c)\geq2r_n\right]\\
&\geq\exp\left\{-\widetilde C_3\left(
(\log(\e+n))^{\frac{1-\varepsilon}{\overline\varepsilon}-\frac{\iota}{1-\delta}}
+(\log(\e+n))^{1-\varepsilon}\right)\right\}
\end{aligned}
\end{equation*}
for some $\widetilde C_3<\infty$ depending only on $\beta$. Moreover, Proposition \ref{uppertail-up} with $\eta=1$ implies
\begin{equation*}
\mathds P\left[\max_{z\in[-a_{n-1/2},a_{n-1/2}]}R(0,z)>\frac{r_n}{\lambda_n}\right]
\leq c_{\eta=1}^{-1}\exp\left\{-c_{\eta=1}\frac{r_n}{\lambda_na_{n-1/2}^\delta}\right\}
\leq\widetilde c_4n^{-2}
\end{equation*}
for all sufficiently large $n$ and some $\widetilde c_4>0$. Applying these estimates to \eqref{prob-H12H21}, we obtain
\begin{equation*}\label{prob-H12H21-2}
\begin{aligned}
&\I_{H_n\cap\widetilde H_n}\mathds P\left[
H_{1,2}(r_n,\lambda_n)\cap H_{2,1}(r_n,h_2(n))\mid\mathcal F_{n-1}\right]\\
&\geq\I_{H_n\cap\widetilde H_n}\left(
\exp\left\{-\widetilde C_3\left(
(\log(\e+n))^{\frac{1-\varepsilon}{\overline\varepsilon}-\frac{\iota}{1-\delta}}
+(\log(\e+n))^{1-\varepsilon}\right)\right\}
-\widetilde c_4n^{-2}\right).
\end{aligned}
\end{equation*}

To estimate the probability of the event $H_{1,1}(r_n)^c$, we choose $\widetilde{\eta}\in (0,\frac{1}{1-\delta})$ such that  
\begin{equation}\label{def-etatilde}
  \widetilde{\eta}>\frac{1}{1-\delta}\left(\bar{\varepsilon}\vee \left(1-\frac{\iota\bar{\varepsilon}}{(1-\varepsilon)(1-\delta)}\right)\right).  
\end{equation}
Note that, based on the definitions of $\overline{\varepsilon}$ and $\iota$ in \eqref{def-iota}, such $\widetilde{\eta}$ can indeed be chosen. 
On $\widetilde H_n$, the resistance from $0$ to any vertex in $[-a_{n-1},a_{n-1}]$, restricted to this interval, is at most $r_n/\lambda_n$. The edges internal to each of the two remaining intervals are independent of $\mathcal F_{n-1}$. Hence, by symmetry, the resistance triangle inequality and Lemma \ref{twopoint-diam} with $\eta=\widetilde\eta$, we have 
\begin{equation*}\label{prob-H11}
\begin{aligned}
\I_{H_n\cap\widetilde H_n}\mathds P[H_{1,1}(r_n)^c\mid\mathcal F_{n-1}]
&\leq 2\I_{H_n\cap\widetilde H_n}\mathds P\left[
\max_{i,j\in[a_{n-1}+1,w_n]}R_{[a_{n-1}+1,w_n]}(i,j)\geq r_n/3\right]\\
&\leq 2C_{\widetilde\eta}\I_{H_n\cap\widetilde H_n}
\exp\left\{-\left(\frac{r_n}{3w_n^\delta}\right)^{\widetilde\eta}\right\}\\
&=2C_{\widetilde\eta}\I_{H_n\cap\widetilde H_n}
\exp\left\{-3^{-\widetilde\eta}
(\log(\e+n))^{\frac{(1-\varepsilon)(1-\delta)\widetilde\eta}{\overline\varepsilon}}\right\}
\end{aligned}
\end{equation*}
for all sufficiently large $n$. In addition, by \eqref{def-iota} and \eqref{def-etatilde}, 
\[
0<\frac{1-\varepsilon}{\overline\varepsilon}-\frac{\iota}{1-\delta}<1,
\qquad
\frac{(1-\varepsilon)(1-\delta)\widetilde\eta}{\overline\varepsilon}
>\max\left\{1-\varepsilon,
\frac{1-\varepsilon}{\overline\varepsilon}-\frac{\iota}{1-\delta}\right\}.
\]
Consequently, combining these estimates with \eqref{prob-HH}, we get that, for all sufficiently large $n$,
\begin{equation*}\label{prob-H1112}
\begin{aligned}
&\I_{H_n\cap\widetilde H_n}\mathds P\left[
H_{1,1}(r_n)\cap H_{1,2}(r_n,\lambda_n)\cap H_{2,1}(r_n,h_2(n))
\mid\mathcal F_{n-1}\right]\\
&\geq\frac12\I_{H_n\cap\widetilde H_n}
\exp\left\{-\widetilde C_3\left(
(\log(\e+n))^{\frac{1-\varepsilon}{\overline\varepsilon}-\frac{\iota}{1-\delta}}
+(\log(\e+n))^{1-\varepsilon}\right)\right\}.
\end{aligned}
\end{equation*}
Since $H_n\cap\widetilde H_n$ occurs eventually almost surely, it follows that
\begin{equation}\label{prob-H1112-1}
\sum_{n\geq1}\I_{H_n\cap\widetilde H_n}\mathds P\left[
H_{1,1}(r_n)\cap H_{1,2}(r_n,\lambda_n)\cap H_{2,1}(r_n,h_2(n))
\mid\mathcal F_{n-1}\right]=\infty\qquad\mathds P\text{-a.s.}
\end{equation}

Since $H_n\cap\widetilde H_n$ occurs for all sufficiently large $n$ almost surely, \eqref{prob-H1112-1} implies
\begin{equation*}\label{prob-H1112-2}
\sum_{n\geq1}\mathds P\left[
H_{1,1}(r_n)\cap H_{1,2}(r_n,\lambda_n)\cap H_{2,1}(r_n,h_2(n))
\mid\mathcal F_{n-1}\right]=\infty\qquad\mathds P\text{-a.s.}
\end{equation*}
Recall that $H_{1,1}(r_n),\ H_{1,2}(r_n,\lambda_n),\ H_{2,1}(r_n,h_2(n))\in \mathcal{F}_n$ for all $n\in \mathds{N}$. Consequently, by \cite[Theorem~4.5.5]{Durrett}, one has that $\mathds{P}$-a.s.\ $\{H_{1,1}(r_n)\cap H_{1,2}(r_n,\lambda_n)\cap H_{2,1}(r_n,h_2(n)) \text{ i.o.}\}$ occurs. Furthermore, by combining this with \eqref{dcp-H1}, we get that $\mathds{P}$-a.s.\ $\{ H'_{1}(r_n,\lambda_n)\cap H_{2,1}(r_n,h_2(n)) \text{ i.o.}\}$ occurs. 
Together with the eventual occurrence of $H_{2,2}(r_n,h_2(n))\cap H_3(r_n,\lambda_n,h_3(n))$, this verifies Assumption \ref{assumpC}(1).
\end{proof}

By using the arguments in the verification of Assumption \ref{assumpC}(1), we can prove Proposition \ref{Vn-lu}(2) immediately.

\begin{proof}[Proof of Proposition \ref{Vn-lu}(2)]
By the preceding verification of Assumption \ref{assumpC}(1), we see that almost surely,
\[
H_2(r_n,h_2(n))\quad\text{occurs infinitely often}.
\]
On this event, $\bm V(0,r_n)\leq r_n^{1/\delta}/h_2(n)$. We now take
\[
\alpha_2=\frac{(1-\varepsilon)(1-\delta)}{\delta}>0,
\]
where $\varepsilon$ is chosen in \eqref{eta}. By \eqref{v1-n} and $r_n=M^{n^2}$, one has 
\[
h_2(n)\asymp(\log n)^{\alpha_2}
\asymp(\log\log r_n)^{\alpha_2}.
\]
Consequently, we obtain that for some $\widetilde c_1<\infty$ depending only on $\beta$,
\[
\liminf_{n\to\infty}
\frac{\bm V(0,r_n)(\log\log r_n)^{\alpha_2}}{r_n^{1/\delta}}
\leq \widetilde c_1\quad\mathds P\text{-a.s.}
\]
Replacing $r_n$ by $\lfloor r_n\rfloor$ and using the monotonicity of $\bm V(0,\cdot)$, we obtain the desired estimate. 
\end{proof}

We now present the 
\begin{proof}[Proof of Theorem \ref{osc-LRP}]
Recall that $\delta\in (0,1)$ is the exponent defined in Theorem \ref{eff-tightness}. From Proposition \ref{Vn-lu}(1), 
Assumption \ref{assumpB}(1) holds with $d_f =1/\delta$ and $h_1(n) =\widetilde{c}_1(\log\log (\e^{\e}+n))^{\alpha_1}$, where $\alpha_1>0$ is the constant as defined in Proposition \ref{Vn-lu} and $\widetilde{c}_1>0$ is the normalizing constant such that $h_1(1)=1$.  
Additionally, Proposition \ref{Vn-lu}(2) implies that Assumption \ref{assumpB}(2) holds with $d_f = 1/\delta$ and $h_2(n) =\widetilde{c}_2(\log\log (\e^{\e}+n))^{\alpha_2}$, where $\alpha_2>0$ is the constant as defined in Proposition~\ref{Vn-lu} and $\widetilde{c}_2>0$ is the normalizing constant such that $h_2(1)=1$.
Finally, by Lemmas \ref{lem-rlow-hklu}, \ref{V-lower} and \ref{up-V}, we get that Assumption \ref{assumpA} hold with $h(\lambda)=\lambda^{-\min\{q_1,q_2\}}$. 
Hence, Theorem \ref{osc-LRP}(1), (2) and (4) follow from Theorems \ref {thm-tight}, \ref{hk-limif} and Proposition~\ref{thm:tightnot}.

By the verification of Assumption \ref{assumpC}(1) above, we have $\Phi_2(n)=G(T^{-1}(n))\asymp(\log\log n)^{\alpha_*}$ for some $\alpha_*>0$. Thus \eqref{limsuppt} follows from Theorem \ref{hk-limsup2}.
Finally, Proposition \ref{prop-assumpC2} verifies Assumption \ref{assumpC}(2), so \eqref{limsupext} follows from Theorem \ref{hk-limsup2}.  
\end{proof}

\section{Proof of 
the results in Section \ref{main-results}}
\label{sect-proofmr}

To prove the results presented in Section \ref{main-results}, we will first introduce some relevant results 
from \cite{BK06} concerning random walks on deterministic weighted graphs (see Subsection \ref{subsect-hk}) and on random graphs (see Subsection \ref{subsect-hkRG}). 
We will then apply these results to complete the proof of the conclusions in Section \ref{main-results}.

\subsection{Heat kernel estimates on graphs}\label{subsect-hk}

Consider $\mathcal{G}=(V,E)$ as a locally finite, connected undirected graph that includes a marked vertex $0\in V$, where  $\{\mu_e\}_{e\in E}$ represent the positive conductances on the graph $\mathcal{G}$. The measure $\mu$ on $V$ is induced by $\{\mu_e\}_{e\in E}$ as defined in \eqref{def-meas}. 
Recall that $X=((X_t)_{t\geq 0}, \mathbf{P}_x,x\in V)$ denotes the continuous-time random walk on $\mathcal{G}$ with generator defined by \eqref{def-generator}, and $p_t(\cdot,\cdot)$ represents the transition density of $X$ with respect to $\mu$.  Additionally, we recall that $B(\cdot,\cdot)$ denotes the ball under the resistance distance as defined in \eqref{def-ball}, and $\bm V(\cdot,\cdot)$ represents the volume of the ball $B(\cdot,\cdot)$ with respect to $\mu$. 

\medskip

In this subsection, we state three basic estimates from \cite{BK06} concerning heat kernel estimates related to the volume $\bm V(\cdot,\cdot)$. 
For readers' convenience, we give proof. 
Note that the results in \cite{BK06}
are stated using the graph distance, but the same proof works for resistance metric. 
Similar results for discrete time Markov chains are given in \cite{BJKS08,Kum14,KM08} as well.
In particular, corresponding results are stated for general metric (not necessarily the graph distance) in \cite{Kum14,KM08}. 

The following proposition provides a upper bound of the heat kernel $p_t(\cdot,\cdot)$. 

\begin{proposition}{\rm (\cite[Theorem 4.1]{BK06})}\label{prop-hkV}
For all $x\in V$ and all $r\geq 1$, we have 
\begin{equation}\label{Gamma-V}
  p_{2t}(x,x)\leq \frac{2}{\bm V(x,r)}\vee \frac{2r}{t},\quad t>0. 
  \end{equation}
Hence, 
  \begin{equation*}
    p_{2r\bm V(x,r)}(x,x)\leq \frac{2}{\bm V(x,r)}.
  \end{equation*}
\end{proposition}
\begin{proof}
  For fixed $x\in V$ and $r\geq 1$, we write $B(r)=B(x,r)$ and $\bm V(r)=\bm V(x,r)$ for simplicity. For any $t>0$, let $f_t(\cdot)=p_t(x,\cdot)$ be a function on $V$ and let $\Gamma(t)=\mu(f^2_t)$. By the symmetry of the heat kernel $p_t(\cdot,\cdot)$ and the Markov property of the random walk $X$, it is clear that 
  \begin{equation*}
    \Gamma(t)=\sum_{y\in V}p_t(x,y)p_t(y,x)\mu_y=p_{2t}(x,x)=f_{2t}(x).
  \end{equation*}
  Since
  \begin{equation*}
    \sum_{y\in B(r)}f_t(y)\mu_y=\sum_{y\in B(r)}\mathbf{P}_x\left[X_t=y\right]\leq 1,
  \end{equation*}
  there exists $y\in B(r)$ (depends on $t$ and $r$) with $f_t(y)\leq 1/\bm V(r)$.
  Then by \cite[Lemma 3.1]{BK06}, we have 
  \begin{equation*}
    \frac{1}{2}f_t(x)^2\leq f_t(y)^2+|f_t(x)-f_t(y)|^2\leq \frac{1}{\bm V(r)^2}+R(x,y)\mathcal{D}(f_t,f_t)\leq \frac{1}{\bm V(r)^2}+r\mathcal{D}(f_t,f_t).
  \end{equation*}
  Combining this with the definition of $\Gamma(t)$, we arrive at
  \begin{equation*}
    \frac{\d\Gamma(t)}{\d t} =-2\mathcal{D}(f_t,f_t)\leq \frac{2\bm V(r)^{-2}-\Gamma(t/2)^2}{r}.
  \end{equation*}
  Note that $-\Gamma(s/2)\leq -\Gamma(t)$ for $t\leq s\leq 2t$. Integrating the above inequality from $t$ to $2t$, we obtain that
  \begin{equation*}
    \Gamma(2t)-\Gamma(t)\leq 2tr^{-1} \bm V(r)^{-2}-tr^{-1}\Gamma(t)^2.
  \end{equation*}
  This implies that
  \begin{equation*}
    t\bm V(r)^2\Gamma(t)^2\leq 2t+r\bm V(r)^2\Gamma(t)\leq (4t)\vee\left(2r\bm V(r)^2\Gamma(t)\right).
  \end{equation*}
  Consequently, we obtain \eqref{Gamma-V}. 
  Hence, the proof is complete by taking $t=r\bm V(r)$.
\end{proof}

We also need the following relation of the heat kernel and the effective resistance. 

\begin{lemma}{\rm (\cite[Lemma 4.3]{BK06})}\label{hk-R}
  For fixed $x,y\in V$ and $t>0$, we have
  \begin{equation*}
    \left|\frac{p_t(x,y)}{p_t(x,x)}-1\right|^2\leq \frac{R(x,y)}{tp_t(x,x)}.
  \end{equation*}
\end{lemma}
\begin{proof}
  For convenience, let $f_t(y)=p_t(x,y)$ and $\psi(t)=\|f_t\|_2^2$. It is clear that $\psi'(t)=-2\mathcal{D}(f_t,f_t)$ and $\mathcal{D}(f_t,f_t)$ is decreasing with $t$. Then we have
  \begin{equation*}
    t\mathcal{D}(f_t,f_t)\leq 2\int_{t/2}^t\mathcal{D}(f_s,f_s)\d s=\psi(t/2)-\psi(t)\leq \psi(t/2)= f_t(x).
  \end{equation*}
  Therefore, by \cite[Lemma 3.1]{BK06}, we get that
  \begin{equation*}
    |f_t(x)-f_t(y)|^2\leq R(x,y)\mathcal{D}(f_t,f_t)\leq R(x,y)f_t(x)/t.  \qedhere
  \end{equation*}
\end{proof}

\begin{proposition}{\rm (\cite[Proposition 4.4]{BK06})}\label{lem-hklu}
  For $x_0\in V$, $r\geq 1$, $\lambda\geq 1$ and $\varepsilon \leq 1/(2\lambda)$, we have 
  \begin{itemize}
  \item[\rm (1)]
  $
    \mathbf{E}_x\left[\tau_{B(x_0,r)}\right]\leq 2r\bm V(x_0,r)$ for all  $x\in B(x_0,r)$.
   \end{itemize}
  Suppose further that
    \begin{equation}\label{cond-R(x,B)}
    R(x,B(x_0,r)^c)\geq r/\lambda\quad \text{for all }x\in B(x_0,\varepsilon r). 
  \end{equation}
    Then we have 
  \begin{itemize}
    \item[\rm (2)] For $x\in B(x_0,\frac{1}{2}\varepsilon r)$,
    \begin{equation*}
      \mathbf{E}_x\left[\tau_{B(x_0,r)}\right]\geq \frac{r}{4\lambda}\bm V(x_0,\varepsilon r/2).
    \end{equation*}
    \item[\rm (3)]For $x\in B(x_0,\frac{1}{2}\varepsilon r)$,
    \begin{equation*}
      \mathbf{P}_x\left[\tau_{B(x_0,r)}\leq t \right]\leq 1-\frac{\bm V(x_0,\varepsilon r/2)}{8\lambda \bm V(x_0,r)}+\frac{t}{2r\bm V(x_0, r)}.
    \end{equation*}

    \item[\rm (4)] There exists a constant $c_1>0$ such that  for all $t<\frac{r}{8\lambda}\bm V(x_0,\varepsilon r/2)$ and $x\in B(x_0,\varepsilon r/2)$,
    \begin{equation*}
      p_{2t}(x,x)\geq \frac{c_1\bm V(x_0,\varepsilon r/2)^2}{\lambda^2\bm V(x_0, r)^3}.
    \end{equation*}
  \end{itemize}
\end{proposition}

\begin{proof}
Throughout the proof, we write $B$, $B'$ for $B(x_0,r)$ and $B(x_0,\frac{1}{2}\varepsilon r)$, respectively.
For convenience, let $X^B$ denote the random walk $X$ killed outside $B$, and let $p^B_t(\cdot,\cdot)$ be the transition density of $X^B_t$ with respect to $\mu$. We also denote $g_B(\cdot,\cdot)$ as the Green kernel of $X^B$, that is, $g_B(x,y)=\int_0^\infty p^B_t(x,y)\d t$. For any $x\in B$, define $q_B^x(y)=g_B(x,y)/g_B(x,x)$. Then it can be checked that $q_B^x$ is the equilibrium potential for $R(x, B^c)$, i.e. 
  \begin{equation}\label{R-g}
    R(x,B^c)^{-1}=\mathcal{D}(q_B^x,q_B^x)=g_B(x,x)^{-1}, 
  \end{equation}
  see \cite[(3.5)]{BK06} for more details. 
  
  (1) For any $x\in B$, since $R(x, B^c)\leq 2r$ by the triangle inequality of resistance, one has
  \begin{equation*}
    \mathbf{E}_x\left[\tau_{B}\right]=\sum_{y\in B}g_B(x,y)\mu_y\leq g_B(x,x)\bm V(x_0,r)=R(x,B^c)\bm V(x_0,r)\leq 2r\bm V(x_0,r).
  \end{equation*}

  (2) By \eqref{R-g} and \cite[Lemma 3.1]{BK06}, we can see that for $x,y\in B'$, 
\begin{equation*}
  |1-q_B^x(y)|^2\leq R(x,y) R(x,B^c)^{-1}\leq \lambda \varepsilon \leq 1/2.
\end{equation*}
This implies that $q_B^x(y)\geq 1-2^{-1/2}\geq 1/4$. Therefore, for $x\in B'$ we have 
  \begin{equation*}
    \mathbf{E}_x\left[\tau_{B}\right]\geq \sum_{y\in B'}g_B(x,y)\mu_y= \sum_{y\in B'}g_B(x,x)q_B^x(y)\mu_y\geq \frac{1}{4}\bm V(x_0,\varepsilon r/2)R(x,B^c)\geq\frac{r}{4\lambda}\bm V(x_0,\varepsilon r/2), 
  \end{equation*}
 where the last inequality is from \eqref{cond-R(x,B)}. 
 
(3) By (1), (2) and the Markov property, we have that for all $x\in B'$ and all $t>0$,
\begin{equation*}
  \frac{r}{4\lambda}\bm V(x_0,\varepsilon r/2)\leq \mathbf{E}_x\left[\tau_{B}\right]\leq t+\mathbf{E}_x\left[\I_{\{\tau_B>t\}}\mathbf{E}_{X_t}[\tau_B]\right]\leq t+2r\bm V(x_0,r)\mathbf{P}_x\left[\tau_B>t\right].
\end{equation*}

(4) By (3) we get that
\begin{equation*}
  \mathbf{P}_x\left[X_t\in B\right] \geq \mathbf{P}_x\left[\tau_B>t\right]\geq \frac{\bm V(x_0,\varepsilon r/2)}{8\lambda \bm V(x_0,r)}-\frac{t}{2r\bm V(x_0, r)}.
\end{equation*}
Consequently, for $t\leq \frac{r}{8\lambda}\bm V(x_0,\varepsilon r/2)$, there exists a constant $\widetilde{c}_1>0$ such that
\begin{equation}\label{Prob-X_t}
  \mathbf{P}_x\left[X_t\in B\right]\geq  \frac{\widetilde{c}_1\bm V(x_0,\varepsilon r/2)}{\lambda \bm V(x_0,r)}.
\end{equation}
Furthermore, by the Chapman-Kolmogorov equation and the Cauchy-Schwarz inequality, we obtain 
\begin{equation*}
  \mathbf{P}_x\left[X_t\in B\right]^2=\left(\sum_{y\in B}p_t(x,y)\mu_y\right)^2\leq \bm V(x_0,r)\sum_{y\in B}p_t(x,y)^2\mu_y\leq p_{2t}(x,x)\bm V(x_0,r).
\end{equation*}
Combining this with \eqref{Prob-X_t} gives the desired result.
\end{proof}

\subsection{Heat kernel estimates on random graphs}\label{subsect-hkRG}
We now consider a random graph $\mathcal{G}(\omega)=(V(\omega), E(\omega)),\ \omega\in \Omega$ under the framework established in Section \ref{main-results}. 
To begin, We introduce the following definition of a good ball in the graph $\mathcal{G}(\omega)$. 
Recall that $d_f>0$ and $q_1,q_2\geq 1$ are constants defined in Assumption \ref{assumpA}. 

\begin{definition}\label{def-goodball}
  For $x\in V$ and $r,\lambda\geq 1$, we say that $B(x,r)$ is $\lambda$-good if 
  \begin{itemize}
  \item[(1)] $R(y,B(x,r)^c)\geq 8r/\lambda$ for all $y\in B(x,r/\lambda)$. 
  
    \item[(2)] $\lambda^{-q_1}r^{d_f}\leq \bm V(x,r)\leq \lambda^{q_2}r^{d_f}$. 
  
    \item[(3)] $\bm V(x,r/\lambda)\geq \lambda^{-(q_1+d_f)}r^{d_f}$. 
    
    \item[(4)] $\bm V(x,r/\lambda^2)\geq \lambda^{-(q_1+2d_f)}r^{d_f}$.  
  \end{itemize}
\end{definition}
Note that $B(x,r)$ that satisfies $(1)', (2), (3), (4)$ is also $\lambda$-good where 
 \begin{itemize}
  \item[(1)'] $R(x,B(x,r)^c)\geq 9r/\lambda$.
  \end{itemize}
 Indeed, by the same reason as in \eqref{eq:resttri}, $(1)'$ implies $(1)$.

Under Assumption \ref{assumpA}, we can show that the ball $B(x,\cdot)$ is good with high probability. 

\begin{lemma}\label{prob-good}
Suppose that $x\in V$ satisfies Assumption \ref{assumpA}. Then there exists a constant $C_1<\infty$ such that for all $r\ge 1$ and all $\lambda\geq 1$, we have 
\begin{equation*}
\mathds{P}\left[B(x,r)\text{ is not $\lambda$-good}\right]\leq  
C_1h(\lambda/14),   
\end{equation*}
where $h$ is the decreasing function defined in Assumption \ref{assumpA}. 
\end{lemma}

\begin{proof}
For fixed $x\in V$ and $\lambda\geq 14$, 
By (A2) in Assumption \ref{assumpA},  
similarly to \eqref{eq:resttri}, we have 
\begin{equation*}\label{prob-1}
\begin{aligned}
  &\mathds{P}\left[R(y,B(x,r)^c)\geq 8r/\lambda\text{ for all }y\in B(x,r/\lambda)\right]\\
  &\geq  \mathds{P}\left[R(x,B(x,r)^c)\geq \frac{7r}{\lambda}\right]
    \geq 1-h(\lambda/14). 
  \end{aligned}
\end{equation*}
Combining this with (A1) in Assumption \ref{assumpA} gives the desired result. 
\end{proof}

In the good ball, the following heat kernel estimates hold, 
which are generalizations of \cite[Theorem 4.6 (b) and (c)]{BK06}.

\begin{lemma}\label{lem-hklu2}
For $x_0\in V$, $r\geq 1$ and $\lambda\geq 1$, suppose that $x_0$ satisfies Assumption \ref{assumpA} and $B(x_0,r)$ is $\lambda$-good. Let $I=I(\lambda,r)=[\lambda^{-2(1+d_f)}r^{1+d_f}, \lambda^{-(1+2d_f)}r^{1+d_f}]$.
\begin{itemize}
  \item[\rm(1)] For each $K\geq 0$,
  \begin{equation*}
    p^\omega_{2t}(x_0,y)\leq (2+\sqrt{2K})\lambda^{q_1+1} t^{-\frac{d_f}{1+d_f}}\quad \text{for all }t\in I,\ y\in B(x_0,Kt^{\frac{1}{1+d_f}}).
  \end{equation*}
  \item[\rm (2)] For $x\in B(x_0, r/\lambda)$, there exist constants $c_2,c_3>0$ such that
  \begin{equation*}
    p^\omega_{2t}(x,y)\geq c_2\lambda^{-(2q_1+2+3q_2)} t^{-\frac{d_f}{1+d_f}}
    \quad \text{if}\quad R(x,y)\leq c_3\lambda^{-(4+2q_1+4d_f+3q_2)}r,\ t\in I.
  \end{equation*}
\end{itemize}
\end{lemma}
\begin{proof}
Throughout the proof, we will omit the superscript of $p_{\cdot}^\omega(\cdot,\cdot)$ for simplicity. The proof is very similar to that of \cite[Theorem 4.6 (b) and (c)]{BK06}. 

  (1) Let $t_1=(r/\lambda^2)^{1+d_f}$. By (4) in Definition \ref{def-goodball} and using \eqref{Gamma-V} with $t=r^{1+d_f}$, we can see that for $t\in I$,
  \begin{equation*}
    p_{2t}(x_0,x_0)\leq p_{2t_1}(x_0,x_0)\leq \frac{2}{\bm V(x_0,t_1^{\frac{1}{1+d_f}})}\vee \frac{2t_1^{\frac{1}{1+d_f}}}{t_1}
    \leq2 \lambda^{q_1} t_1^{-\frac{d_f}{1+d_f}}\leq 2 \lambda^{q_1+\frac{d_f}{1+d_f}} t^{-\frac{d_f}{1+d_f}}\leq 2 \lambda^{q_1+1} t^{-\frac{d_f}{1+d_f}}. 
  \end{equation*}
  Combining this with Lemma \ref{hk-R}, we find that for $t\in I$ and $y\in B(x_0,Kt^{\frac{1}{1+d_f}})$,
  \begin{equation*}
  \begin{aligned}
    p_{2t}(x_0,y)&\leq p_{2t}(x_0,x_0)+|p_{2t}(x_0,x_0)-p_{2t}(x_0,y)|\\
    &\leq p_{2t}(x_0,x_0)+ \sqrt{R(x_0,y)p_{2t}(x_0,x_0)/t}\\
    &\leq  p_{2t}(x_0,x_0)+\sqrt{Kt^{-\frac{d_f}{1+d_f}}p_{2t}(x_0,x_0)}\leq (2+\sqrt{2K})\lambda^{q_1+1} t^{-\frac{d_f}{1+d_f}}.
    \end{aligned}
  \end{equation*}

  (2) Let $\widetilde{\lambda}=\lambda/14$ and $\varepsilon=4/\lambda$. By Definition \ref{def-goodball}(3), we can see that
  \begin{equation*}
    \bm V(x_0,\varepsilon r/2)\geq \bm V(x_0,r/\lambda)\geq \lambda^{-(q_1+d_f)}r^{d_f}.
  \end{equation*}
From this, $q_1\leq d_f$ and Proposition \ref{lem-hklu}(4) with $\lambda=\widetilde{\lambda}$, we obtain that for $x\in B(x_0,r/\lambda)$ and $t\in I$,
\begin{equation*}
\begin{aligned}
  p_{2t}(x,x)&\geq \frac{\widetilde{c}_1{\bm V}(x_0,\varepsilon r/2)^2}{\lambda^2{\bm V}(x_0, r)^3}\geq \widetilde{c}_1\left(\lambda^{-(q_1+d_f)}r^{d_f}\right)^2\cdot \lambda^{-2}\left(\lambda^{q_2}r^{d_f}\right)^{-3}\\
  &=\widetilde{c}_1\lambda^{-(2q_1+2d_f+2+3q_2)}r^{-d_f}
  \geq\widetilde{c}_1\lambda^{-(2q_1+4d_f+2+3q_2)} t^{-\frac{d_f}{1+d_f}}.
  \end{aligned}
\end{equation*}
Thus, by Lemma \ref{hk-R}, for $R(x,y)\leq \widetilde{c}_1\lambda^{-(4+2q_1+4d_f+3q_2)}r$, 
we have
  \begin{equation*}
    \left|\frac{p_{2t}(x,y)}{p_{2t}(x,x)}-1\right|^2\leq \frac{R(x,y)}{2tp_{2t}(x,x)}\leq \frac{R(x,y)\lambda^{4+2q_1+4d_f+3q_2}}{2\widetilde{c}_1 r} 
     \leq \frac{1}{2}.
  \end{equation*}
  Hence, the proof is complete by some straightforward calculations. 
\end{proof}

\subsection{Proof of Theorems \ref{thm-tight}, \ref{hk-limif}, \ref{hk-limsup2}
and Proposition \ref{thm:tightnot}. 
} \label{proof-mr}

We start with the 

\begin{proof}[Proof of Theorems \ref{hk-limif}]
(1) Recall that $d_f>0$ is the constant defined in Assumption \ref{assumpB}(1).  
For fixed $n\geq 1$, 
let $a_n$ be the random number such that $\bm V(0,r_n)=a_nr_n^{d_f}$  
where $\{r_n\}_{n\geq 1}$ is a sequence in Assumption \ref{assumpB}(1), 
and let $t_n=r_n\bm V(0,r_n)=a_nr_n^{1+d_f}$. Then by Proposition \ref{prop-hkV}, we get that for any $\omega\in \Omega$, 
\begin{equation}\label{hk-est1}
  p_{2t_n}^\omega(0,0)\leq \frac{2}{\bm V(0,r_n)(\omega)}=\frac{2}{a_nr_n^{d_f}}=2a_n^{-\frac{1}{1+d_f}}t_n^{-\frac{d_f}{1+d_f}}.
\end{equation}
Moreover, from Assumption \ref{assumpB}(1), we can see that $\mathds{P}$-a.s.\ $a_n\geq h_1(r_n)$ for infinitely many $n$. Define $\Psi(x):=xh_1^{-1}(x)^{1+d_f}$, where $h_1^{-1}$ is 
the inverse function of $h_1$. Then $a_n\geq h_1(r_n)$ implies $\Psi(a_n)=
a_nh_1^{-1}(a_n)^{1+d_f}\ge a_nr_n^{1+d_f}=t_n$. Now define $\Phi_1(x):=\Psi^{-1}(x)^{1/(1+d_f)}$, 
then $a_n\ge \Phi_1(t_n)^{1+d_f}$, hence $a_n^{-\frac{1}{1+d_f}}\le \Phi_1(t_n)^{-1}$. 
It is easy to check $\Phi_1\in {\bf SV}$.
Applying this to \eqref{hk-est1}, we obtain that $\mathds{P}$-a.s.\ 
\begin{equation*}
  p_{2t_n}^\omega(0,0)\leq 2\Phi_1(t_n)^{-1}t_n^{-\frac{d_f}{1+d_f}}
\end{equation*}
for infinitely many $n$. Hence, the proof is complete. \\
(2) Let $\{r_n\}_{n\geq 1}$ be a sequence in Assumption \ref{assumpB}(2) and 
$a_n$ be the random number such that $\bm V(0,r_n)=a_nr_n^{d_f}$ as before. Then by Proposition \ref{lem-hklu}(1) and Assumption \ref{assumpB}(2), it holds that $\mathds{P}$-a.s.\ 
\[
 \mathbf{E}_{x_0}\left[\tau_{B(x_0,r_n)}\right]\leq 2r_n\bm V(x_0,r_n)
 \le 2c_2r_n^{	1+d_f} /h_2(r_n)
\]
for infinitely many $n$. Hence, the proof is complete. 
\end{proof}

Using Lemmas \ref{prob-good} and \ref{lem-hklu2}, we can present the 

\begin{proof}[Proof of Theorem \ref{thm-tight}]
For fixed $x\in V$, $r>0$ and $\lambda\geq 1$, from Lemma \ref{lem-hklu2}, we can see that if $B(x,r)$ is $\lambda$-good, there exist constants $\widetilde{c}_1,\widetilde{c}_2>0$ such that for $t=\lambda^{-2(1+d_f)}r^{1+d_f}$, 
\begin{equation*}
  \widetilde{c}_1\lambda^{-(2q_1+2+3q_2)}\leq t^{\frac{d_f}{1+d_f}}p^\omega_{2t}(x,x)\leq \widetilde{c}_2\lambda^{q_1+1}. 
\end{equation*}
Also, by Proposition \ref{lem-hklu}, 
we can see that if $B(x,r)$ is $\lambda$-good and $\lambda\ge 4$, then 
\begin{equation*}
  \frac 14 r^{1+d_f}\lambda^{-(1+q_1+2d_f)}\le 
  \mathbf{E}_x\left[\tau_{B(x,r)}\right]\le 2r^{1+d_f}\lambda^{q_2}.
\end{equation*}
Consequently, for any sufficiently small $\varepsilon\in (0,1)$, 
we take $\lambda \geq 4$ 
to be the smallest integer that satisfies 
\begin{equation*}
  \varepsilon \leq 
   (\frac 14 \lambda^{-(1+q_1+2d_f)})\wedge
  (\widetilde{c}_1\lambda^{-(2q_1+2+3q_2)})<
  (2\lambda^{q_2})\vee(\widetilde{c}_2\lambda^{q_1+1})
  \leq \varepsilon^{-1}.
\end{equation*}
Then by Lemma \ref{prob-good}, we get that 
\begin{eqnarray*}
\mathds{P}\left[\varepsilon \leq
  r^{-\frac{1+\delta}{\delta}}\mathbf{E}_x\left[\tau_{B(x,r)}\right]\leq
  \varepsilon^{-1}\right]\wedge    
\mathds{P}\left[\varepsilon \leq  t^{\frac{d_f}{1+d_f}}p^\omega_{2t}(x,x)\leq \varepsilon^{-1}\right]&\geq& \mathds{P}\left[B(x,r)\text{ is $\lambda$-good}\right]\\
&\geq& 1-C_1h(\varepsilon^{-\eta}), 
\end{eqnarray*}
  where $\eta>0$ is a constant depending only on $q_1$ and $q_2$. 
\end{proof}

We now give the 

\begin{proof}[Proof of Theorem \ref{hk-limsup2}]
(1) By Assumption \ref{assumpC}(1) and Proposition \ref{lem-hklu}(4) with $\varepsilon:= (2\lambda)^{-1}$,  we have   \begin{equation*}
      p_{2t}(0,0)\geq \frac{c_1\bm V(0, r/(4\lambda))^2}{\lambda^2\bm V(0, r)^3}
    \end{equation*}
for all $t<\frac{r}{8\lambda}\bm V(0, r/(4\lambda))$, where $c_1>0$ is the constant defined in Proposition \ref{lem-hklu}(4). 
Now set $T_n=\Big(\frac{r_n}{\lambda_n}\Big)^{1+d_f}h_3(n)^{-1}$ 
as in \eqref{eq:T_ndesu}. Then, 
on $H_1(r_n,\lambda_n)\cap H_2(r_n,h_2(n))\cap H_3(r_n,\lambda_n,h_3(n))$, we have 
$T_n \le \frac{r_n}{c_4\lambda_n}\bm V(0, r_n/(4\lambda_n))$,
where $c_4>0$ is the constant defined in Definition \ref{def-I-1}. 
Hence, we get 
\begin{equation*}
     T_n^{d_f/(1+d_f)} p_{2T_n}(0,0)\geq  T_n^{d_f/(1+d_f)}\frac{c_1
     \left(c_4\Big(\frac{r_n}{\lambda_n}\Big)^{d_f}h_3(n)^{-1}\right)^2}{\lambda_n^2\left(c_3r_n^{d_f}h_2(n)^{-1}\right)^3}
\ge 
\widetilde{C}\frac{h_2(n)^3}{\lambda_n^{2+3d_f}h_3(n)^{2+d_f/(1+d_f)}}=\widetilde{C}G(n),  
    \end{equation*}
  for some $\widetilde{C}>0$ where 
  $G(n)$ is the function defined in \eqref{condC}.  Since $G(n)= \Phi_2(T_n)$, 
we obtain the desired estimate. In addition, note that by the definitions of $G$ and $T$ in \eqref{condC} and \eqref{eq:T_ndesu} respectively,  
it is easy to check $\Phi_2\in {\bf SV}$. \\
(2) By Assumption \ref{assumpC}(2) and Proposition \ref{lem-hklu}(2) with $\varepsilon:= (2\lambda)^{-1}$,  we have   
\begin{equation*}
   r_n^{-(1+d_f)} \mathbf{E}_{x_0}\left[\tau_{B(x_0,r_n)}\right]\geq  r_n^{-(1+d_f)}\frac{r_n}{4\lambda_n}\bm V(x_0,\frac{r_n}{4\lambda_n})
  \geq \widetilde{C'}\frac{h_3(n)}{\lambda_n^{1+d_f}}
  =\widetilde{C'}\widehat G(n).
  \end{equation*}
  for some $\widetilde{C'}>0$ where 
  $\widehat G(n)$ is the function defined in \eqref{condD2}.  Since 
  $\widehat G(n)= \Phi_3(r_n)$, 
we obtain the desired estimate. It is easy to check $\Phi_3\in {\bf SV}$.
\end{proof}

Finally, we give 

\begin{proof}[Proof of Proposition \ref{thm:tightnot}]
We follow the proof of \cite[Theorem 1.6]{BK06}. 
For any $n\geq 1$, let $H_n= \sup_{0\le s \le 1}  Z^{(n)}_s$. Then, by Proposition \ref{lem-hklu}(1) we arrive at 
\begin{eqnarray*}
\mathbf{P}_0\left[ H_n \le \lambda\right] &=& 
  \mathbf{P}_0 \left[ \sup_{t\le n} R(0,X_t) \le \lambda n^{1/(d_f+1)}\right] 
 = \mathbf{P}_0\left[ \tau_{B(0,\lambda n^{1/(d_f+1)})}  \ge n \right] \\
 &\le & \frac{\mathbf{E}_x[\tau_{B(0,\lambda n^{1/(d_f+1)})}]}n
\le \frac{ 2 \lambda n^{1/(d_f+1)} \bm V(0,\lambda n^{1/(d_f+1)})}{n}.
\end{eqnarray*}
Hence by \eqref{eq:liminfvol}, for any $\lambda>0$, we have that  
$\liminf_{n\to\infty}  \mathbf{P}^{\omega}_0[ H_n \le \lambda] =0$. This shows that
the random variable $H_n$ (and hence $Z^{(n)}$) are not tight.
\end{proof}

\bigskip

\noindent{\bf Acknowledgement.} \rm
We thank David Croydon and Haojie Hou 
for helpful discussions. 
L.-J.\ Huang is partially supported by National Key R\&D Program of China No.\ 2023YFA1010400, the National Natural Science Foundation of China No.\ 12622110, 12471136, and the Fujian Natural Science Foundation Youth Project No.\ 2026J0141. 
T.\ Kumagai is partially supported by JSPS KAKENHI Grant Numbers
22H00099, 23KK0050 and JST Moonshot R\&D Grant Number JPMJMS25A5.

\noindent{\bf Declaration of AI Use.} \rm
The authors acknowledge the use of AI tools in assisting with the proof of Proposition \ref{prop-assumpC2}, as well as in language polishing, and checking the manuscript. 
All mathematical arguments and results have been carefully reviewed and verified by the authors, who take full responsibility for the accuracy and integrity of the manuscript.

\bibliographystyle{plain}
\bibliography{resistance-ref2}

\end{document}